\documentclass[11pt,leqno,fleqn]{article}
\usepackage{epsf,amsmath,amssymb,amscd,hyperref,srcltx}
\usepackage{amsfonts}
\usepackage{counters}
\usepackage{psfrag}
\newcounter{figuren}

\newcommand{\dy}[2]{%
\refstepcounter{equation}%
\LABEL{#1}%
\begin{list}{}{
\topsep 5mm
\leftmargin 18mm
\rightmargin 0cm
\itemsep 0mm
\listparindent 0mm
\parsep 0mm
\itemsep 0mm
\labelsep 0mm
\labelwidth 18mm
}%
\item[\rm (\theequation)\hfill]
#2
\end{list}%
}
\newcommand{\dyz}[1]{%
\refstepcounter{equation}%
\begin{list}{}{
\topsep 5mm
\leftmargin 18mm
\rightmargin 0cm
\itemsep 0mm
\listparindent 0mm
\parsep 0mm
\itemsep 0mm
\labelsep 0mm
\labelwidth 18mm
}%
\item[\rm (\theequation)\hfill]
#1
\end{list}%
}
\newcommand{\de}[2]{\dy{#1}{\raggedright$\displaystyle #2 $}}
\newcommand{\dez}[1]{\dyz{\raggedright$\displaystyle #1 $}}
\newcommand{\leeg}[1]{}

\newcommand{\di}[2]{%
\refstepcounter{equation}%
\LABEL{#1}%
\begin{list}{}{
\topsep 5mm
\leftmargin 10mm
\rightmargin 0cm
\itemsep 0mm
\listparindent 0mm
\parsep 0mm
\labelsep 1mm
\labelwidth 10mm
}%
\item[\rm (\theequation)\hfill]
\begin{list}{}{
\topsep 0mm
\leftmargin 8mm
\rightmargin 0mm
\itemsep 0mm
\listparindent 0mm
\parsep 0mm
\labelsep 1.5mm
\labelwidth 6.5mm
}
#2
\end{list}%
\end{list}%
}
\newcommand{\dt}[3]{\dy{#1}{\begin{tabular}[t]{#2}#3\end{tabular}}}

\newcommand{\nr}[1]{\item[{\rm (#1)}]}
\newcommand{\nrs}[1]{\item[{\rm (#1)}]\vspace{-\itemsep}}
\newcommand{\items}[1]{\item[#1]\vspace{-\itemsep}}
\newcounter{stelling}
\newcommand{\thm}[2]{\setcounter{gevolg}{0}\setcounter{claim}{0}\refstepcounter{stelling}\vspace{4mm}\noindent{\bf Theorem \thestelling.}\label{#1}{\it #2}}

\newalphacounterwithin{gevolg}{stelling}
\newcommand{\cor}[2]{\refstepcounter{gevolg}\setcounter{claim}{0}\vspace{4mm}\noindent{\bf Corollary \thegevolg.}\label{#1}{\it #2}}

\newcounter{hulpstelling}

\newcounter{bewering}
\newcommand{\prop}[2]{\refstepcounter{bewering}\vspace{4mm}\noindent{\bf Proposition \thebewering.}\label{#1}{\it #2}}

\newcounter{claim}

\newalphacounterwithin{subclaim}{claim}

\newcounter{opmerking}
\newcommand{\opm}[1]{\refstepcounter{opmerking}\vspace{4mm}\noindent{\bf Remark \theopmerking.} \label{#1}}

\newcounter{hoofdstuk}

\newcounter{sectie}
\newcounter{subsectie}
\newcounterwithin{ex}{sectie}
\newcommand{\sect}[2]{\refstepcounter{sectie}\setcounter{subsectie}{0}\setcounter{ex}{0}
\section*{\boldmath \thesectie. #2}%
\label{#1}}
\newcommand{\sectz}[1]{\refstepcounter{sectie}\setcounter{subsectie}{0}\setcounter{ex}{0}
\section*{\boldmath \thesectie. #1}%
}
\newcommand{\subsect}[2]{\refstepcounter{subsectie}
\subsection*{\boldmath \thesectie.\thesubsectie. #2}%
\label{#1}}

\newcounter{lit}

\newcommand{\pf}{\vspace{3mm}\noindent{\bf Proof.}\ }

\newcommand{\bx}{\hspace*{\fill} \hbox{\hskip 1pt \vrule width 4pt height 8pt depth 1.5pt \hskip 1pt}

\addvspace{4mm}}

\newcommand{\rf}[1]{{\rm (\ref{#1})}}

\newcommand{\vv}[1]{\vspace{#1}}

\newcommand{\vn}[1]{\vspace{#1} \noindent}

\newcommand{\dist}{\text{\rm dist}}

\newcommand{\CC}{{\cal C}}

\newcommand{\FF}{{\cal F}}

\newcommand{\PP}{{\cal P}}

\newcommand{\RR}{{\cal R}}

\newcommand{\UU}{{\cal U}}

\newcommand{\kint}[2]{\mbox{$\int$}}

\newcommand{\NIET}[1]{}

\newcommand{\LABEL}[1]{\label{#1}}

\newcommand{\oZ}{{\mathbb{Z}}}
\newcommand{\clos}{\mbox{\rm cl}}
\newcommand{\fff}{\mbox{\sl first}}
\begin{document}
\begin{center}
\baselineskip 7mm
{\LARGE\bf Free partially commutative groups, cohomology, and
paths and circuits in directed graphs on surfaces

}

\end{center}
\vspace{1mm}
\begin{center}
{\large
Alexander Schrijver\footnote{ CWI and University of Amsterdam.
Mailing address: CWI, Science Park 123, 1098 XG Amsterdam,
The Netherlands.
Email: lex@cwi.nl.}}

\end{center}

\noindent
{\small {\bf Abstract.}
We show that for each fixed $k$, the problem of finding $k$ pairwise
vertex-disjoint directed paths between given source-sink pairs in a
planar directed graph is solvable in polynomial time.
In fact, it suffices to fix the number of faces
needed to cover all sources and sinks.
Moreover, the method can be extended to any fixed compact orientable
surface (instead of the plane) and to rooted trees (instead of paths).

Our approach is algebraic and is based on cohomology over graph (nonabelian) groups.
More precisely, let $D=(V,A)$ be a directed graph and let $(G,\cdot)$ be a group.
Call two function $\phi,\psi:A\to G$ {\em cohomologous} if there exists a function $p:V\to G$ such that $p(u)\cdot\phi(a)\cdot p(w)^{-1}=\psi(a)$ for each arc $a=(u,w)$.
Now given a function $\phi:A\to G$ we want to find a function
$\psi$ cohomologous to $\phi$ such that each $\psi(a)$ belongs to a
prescribed subset $H(a)$ of $G$.
We give a polynomial-time algorithm for this problem in case $G$
is a graph group and each $H(a)$ is closed
(i.e., if word $xyz$ belongs to $H(a)$ then also word $y$ belongs to $H(a)$).

The method also implies that such a $\psi$ exists, if and only if for each
$s\in V$ and each pair $P,Q$ of (undirected) $s-s$ paths there exists
an $x\in G$ such that $x\cdot\phi(P)\cdot x^{-1}\in H(P)$ and
$x\cdot\phi(Q)\cdot x^{-1}\in H(P)$.
(Here $\phi(P)$ is the product of the $\phi(a)$ over the arcs in $P$.
Similarly, $H(P)$ is the (group subset) product of the $H(a)$.)

}

\sect{s1}{Introduction}

In this paper we show that the following problem, the $k$ {\em disjoint paths problem for directed planar graphs}, is solvable in polynomial time, for any fixed $k$:
\di{d1}{
\item[{\rm given:}]  a planar directed graph $D=(V,E)$ and $k$ pairs
$(r_1,s_1),\ldots,(r_k,s_k)$ of vertices of $D$; 
\items{\rm find:} $k$ pairwise vertex-disjoint directed paths
$P_1,\ldots,P_k$ in $D$, where $P_i$ runs from $r_i$ to $s_i$
($i=1,\ldots,k$).
}

The problem is NP-complete if we do not fix $k$ (even in the
undirected case; Lynch [6]).
Moreover, it is NP-complete for $k=2$ if we delete the planarity
condition (Fortune, Hopcroft, and Wyllie [5]).
This is in contrast to the undirected case (for those believing
NP$\neq$P), where
Robertson and Seymour [10] showed that, for any fixed $k$, the $k$ disjoint paths 
problem is polynomial-time solvable for any graph (not necessarily
planar).

Our algorithm is a `brute force' polynomial-time algorithm.
We did not aim at obtaining the best possible running time bound, as
we presume that there are much faster (but possibly more complicated)
methods for \rf{d1} than the one we describe in this paper.
In fact, recently Reed, Robertson, Schrijver, and Seymour [9]
showed that for undirected planar graphs the $k$ disjoint paths problem
can be solved in {\em linear} time, for any fixed $k$.
This algorithm makes use of methods from Robertson and Seymour's
theory of graph minors.
A similar algorithm for directed planar graphs might exist, but
probably would require extending parts of graph minors theory to the
directed case.

Our method is based on cohomology over free (nonabelian) groups.
For the $k$ disjoint paths problem we use free groups with $k$ generators.
Cohomology is in a sense dual to homology, and can be defined in any directed graph, also if it is not embedded on a surface.
We apply cohomology to an {\em extension} of the planar graph dual of $D$---just using homology to $D$ itself seems not powerful enough.

This approach allows application of the algorithm where the embedding
of the graph in the plane is given in an implicit way, viz.\ by a list of the cycles that bound the faces of the graph.

It also allows a more general application than \rf{d1}.
It is not necessary to fix the number $k$ of pairs $(r_i,s_i)$
but it suffices to fix the number $p$ of faces of $D$ such that each of
$r_1,s_1,\ldots,r_k,s_k$ is incident with at least one of these faces.
The planarity condition can be relaxed to being embeddable in some fixed compact orientable surface.
(A compact orientable surface is any space obtained from the sphere by adding a finite number of `handles'.)
We can restrict for each arc $a$ the connections that can be made over $a$.
Moreover, the method extends to finding rooted trees instead of directed paths.

That is, for any fixed compact orientable surface $S$ and any fixed $p$ we give a polynomial-time algorithm for the following problem:
\di{d2}{
\item[{\rm given:}] a directed graph $D=(V,A)$ embedded on $S$, subsets $A_1,\ldots,A_k$ of $A$, pairs
$(r_1,S_1),\ldots,(r_k,S_k)$, where $r_i\in V$ and $S_i\subseteq V$ ($i=1,\ldots,k$), such that there exist at most $p$ faces such that each vertex in $\{r_1,\ldots,r_k\}\cup S_1\cup\cdots\cup S_k$ is incident with at least one of these faces; 
\items{\rm find:} $k$ pairwise vertex-disjoint rooted trees
$T_1,\ldots,T_k$, where $T_i$ is rooted in $r_i$, covers $S_i$ and contains arcs only in $A_i$ ($i=1,\ldots,k$).
}
There are several other variants of this problem where the methods below are applicable.
We did however not see if our methods extend to compact {\em non}orientable surfaces.

\sect{s1a}{Directed graphs and surfaces}

We give some notation and terminology on directed graphs and surfaces.
Directed graphs may have loops and parallel arcs. 
Nevertheless we sometimes write $a=(u,w)$, meaning that $a$ is an arc from $u$ to $w$.
For each arc $a$ from $u$ to $w$, we define $a^{-1}$ as the reverse arc from $w$ to $u$. (This need not be an arc of $D$ again.)

An ({\em undirected}) {\em path} is a word
\de{a1}{
P=a_1a_2\cdots a_m,
}
where  $a_i=a$ or $a_i=a^{-1}$ for some arc $a$ ($i=1,\ldots,m$), such that the head of $a_{i}$ is equal to the tail of $a_i$ ($i=1,\ldots,m-1$).
We allow $a_{i+1}=a_i^{-1}$ in \rf{a1}.
Moreover we allow the empty path $\emptyset$, where $m=0$.

We call $P$ an $s-t$ path if $s$ is the tail of $a_1$ and $t$ is the
head of $a_m$.
If $s=t$ we call $P$ a {\em cycle}.
If $P$ is as in \rf{a1} then $P^{-1}:=a_m^{-1}\cdots a_1^{-1}$.

\sect{13me03a}{Graph groups}

\subsect{sN1}{Graph groups}

Our method uses the framework of combinatorial group theory, viz. groups defined by generators and relations.
For background literature on combinatorial group theory we refer to 
Magnus, Karrass, and Solitar [8] and Lyndon and Schupp [7];
however, our treatment below is self-contained.

We first give some standard terminology.
Let $g_1,\ldots,g_k$ be `generators'.
Call the elements $g_1,g_1^{-1},\ldots,g_k,g_k^{-1}$ {\em symbols}.
Define $(g_i^{-1})^{-1}:=g_i$.
A {\em word} (of {\em size} $t$) is a sequence $a_1\cdots a_t$ where each $a_j$ is a symbol.
The empty word (of size 0) is denoted by $\emptyset$.
Define $(a_1\cdots a_t)^{-1}:=a_t^{-1}\cdots a_1^{-1}$.

Word $y$ is a {\em segment} of word $w$ if $w=xyz$ for words $x,z$; $y$ is a {\em beginning segment} if $x=\emptyset$, and an {\em end segment} if $z=\emptyset$.
Word $y= a_1\cdots a_t$ is a {\em subword} of word $w$ if $w=x_0 a_1x_1\cdots x_{t-1} a_tx_t$ for some words $x_0,\ldots,x_t$.
It is a {\em proper} subword if $y\neq w$.

Let $g_1,\ldots,g_k$ be generators, and let $E$ be a set of unordered pairs $\{i,j\}$ from $\{1,\ldots,k\}$ with $i\neq j$.
Then the group $G=G_E$ is generated by the generators $g_1,\ldots,g_k$, with relations
\dy{d5a}{
$g_ig_j=g_jg_i$ for each pair $\{i,j\}\in E$.
}
Such a group is called a {\em free partially commutative group} or a {\em graph group}.
(These groups are studied {\em inter alia} in [1],
[4], [13], [18].
However, in this paper we do not use the results of these papers.)

To describe $G$,
call symbols $a$ and $b$ {\em independent} if $a\in\{g_i,g_i^{-1}\}$ and $b\in\{g_j,g_j^{-1}\}$ for some $\{i,j\}\in E$ with $i\neq j$.
So if $a$ and $b$ are independent then $ab=ba$ and $b\neq a^{\pm 1}$.
(It follows from Proposition \ref{pN5b} below that also the converse implication holds.)

By definition, $G$ consists of all words, identifying any two words $w$ and $w'$ if $w'$ arises from $w$ by iteratively:
\di{dN5c}{
\nr{i} replacing $xaa^{-1}y$ by $xy$ or vice versa, where $a$ is a symbol;
\nrs{ii} replacing $xaby$ by $xbay$ where $a$ and $b$ are independent symbols.
}
By {\em commuting} we will mean applying (ii) iteratively.

Note that if $E=\emptyset$  the group $G_E$ is the {\em free group} generated by $g_1,\ldots,g_k$.
If $E$ consists of {\em all} pairs, then $G_E$ is isomorphic to $\oZ^k$.
Let $1$ denote the unit element of $G$.
So $1=\emptyset$.

A {\em perfect matching on} $\{1,\ldots,t\}$ is a partition of
$\{1,\ldots,t\}$ into pairs.
Pairs $\{i,j\}$ and $\{i',j'\}$ are said to {\em cross} if $i<i'<j<j'$ or $i'<i<j'<j$ (assuming without loss of generality $i<j,i'<j'$).

\prop{pN5z}{
For any word $w=a_1\cdots a_t$ one has: $w=1$ if and only if there exists a perfect matching $M$ on $\{1,\ldots,t\}$ such that 
\di{d5f}{
\nr{i} if $\{i,j\}\in M$ then $a_j=a_i^{-1}$;
\nrs{ii} if two pairs $\{i,j\},\{i',j'\}$ in $M$ cross then $a_i$ and $a_{i'}$ are independent.
}
}

\pf
{\em Necessity.}
If $w=\emptyset$ we can take $M=\emptyset$.
Moreover, one easily shows that the existence of $M$ is maintained under the operations \rf{dN5c}.

{\em Sufficiency.}
Let $M$ satisfy \rf{d5f}.
If $w\not\equiv\emptyset$, choose $\{i,j\}\in M$ with $i<j$ and $j-i$ as small as possible.
Then $a_i$ and $a_{i+1}\cdots a_{j-1}$ are independent,
since each of the pairs containing one of $i+1,\ldots,j-1$ should cross $\{i,j\}$.

Hence $w=a_1\cdots a_{i-1}a_{i+1}\cdots a_{j-1}a_{j+1}\cdots a_t$.
Since $M\setminus\{\{i,j\}\}$ directly gives a perfect matching for the right-hand word, we obtain inductively that $w=1$.
\bx
 
We call a word $w$ {\em reduced} if it is not equal (as a word) to $xaya^{-1}z$ for some symbol $a$ independent of $y$.
We say that a symbol $\alpha$ {\em occurs in} an element $x$ of $G$
if $\alpha$ occurs in any reduced word representing $x$.
We say that two words $x$ and $y$ are {\em independent} if any symbol in $x$ and any symbol in $y$ are independent.
(In particular, $b\neq a^{\pm 1}$ for any symbols $a$ in $x$ and $b$ in $y$.)
Note that reducedness is invariant under commuting.

Proposition \ref{pN5z} directly implies:

\prop{pN5a}{
If $w$ is a reduced word and $w=1$, then $w\equiv\emptyset$.}

\pf
If $w\not\equiv\emptyset$ and $w=1$ one shows, as in the proof of Proposition \ref{pN5z}, that $w$ is not reduced.
\bx

So testing if $w=1$ is easy: just replace (iteratively) any segment
$aya^{-1}$ by $y$ where $a$ is a symbol and $y$ is a word independent of $a$.
The final word is empty if and only if $w=1$.
This gives a test for equivalence of words $w$ and $x$: just test if $wx^{-1}=1$. 
So the `word problem' for free partially commutative groups is easy.
(In fact it can be solved in linear time --- see Wrathall [18].)

Proposition \ref{pN5a} also implies the stronger statement:

\prop{pN5b}{
Let $w$ and $x$ be reduced words with $w=x$.
Then word $x$ can be obtained from $w$ by a series of commutings.
}

\pf
We may assume $w\not\equiv\emptyset\not\equiv x$.
Since $wx^{-1}=1$, $wx^{-1}$ is not reduced.
So we can write $w=w'aw''$ and $z=z'az''$ for some symbol $a$ independent
of $w''$ and $z''$.
By commuting we may assume $w''\equiv z''\equiv\emptyset$.
Then $w'$ and $x'$ are reduced equivalent words, and by induction
$w'$ and $x'$ can be obtained from each other by a series of commutings.
\bx

Proposition \ref{pN5b} implies:
\dy{d7a}{
if $xyz=xy'z$ are reduced words then $y'$ can be obtained
from $y$ by commuting
}
(since $y$ and $y'$ are reduced and $y=y'$.)

In particular, all equivalent reduced words have the same size.
So we can define the {\em size} $|x|$ of an element
$x$ in $G$ as the size of any reduced word $w=x$.
Trivially, $|x^{-1}|=|x|$ and $|xy|\leq |x|+|y|$.
Hence the function $\dist(x,y):=|x^{-1}y|$ is a distance function.
Note that $\dist(zx,zy)=\dist(x,y)$ for all $x,y,z$.

We write
\dyz{
$x|y\Longleftrightarrow |xy|=|x|+|y|$.
}
So
\dyz{
$x|y \Longleftrightarrow$
if $x'$ and $y'$ are reduced words
representing $x$ and $y$, then $x'y'$ is a reduced
word representing $xy$.
}
By extension we write:
\de{d12A}{
x_1|x_2|\cdots|x_n \Longleftrightarrow
|x_1x_2\cdots x_n|=|x_1|+|x_2|\cdots+|x_n|.
}

\subsect{s12A}{The partial order $\leq$}

Let $x$ and $y$ be two reduced words.
We write $x\leq y$ if there are reduced words $x'=x$ and $y'=y$ such that $x'$ is a beginning segment of $y'$.
So $x\leq y$ if and only if $x|x^{-1}y$.

Proposition \ref{pN5b} gives:
\dy{d5h}{
if $x$ and $y$ are reduced words such that $x\leq y$ then $y$ can be commuted to $y'$ such that $x$ is a beginning segment of $y'$.
}
This implies: 

\prop{pN5c}{
$\leq$ is a partial order on $G$.
}

\pf
Clearly $x\leq x$ for each $x\in G$, so $\leq$ is reflexive.
To see that $\leq$ is anti-symmetric, let $x\leq y$ and $y\leq x$.
We may assume that $x$ and $y$ are reduced words.
Then \rf{d5h} implies that $x$ and $y$ can be commuted to each other.
So $x=y$.

To see that $\leq$ is transitive, let $x\leq y$ and $y\leq z$, where $x, y$ and $z$ are reduced words.
By \rf{d5h} $z$ can be commuted to $z'$ such that $y$ is beginning segment of $z'$, and $y$ can be commuted to $y'$ such that $x$ is beginning segment of $y'$.
Hence $z'$ can be commuted to $z''$ such that $x$ is beginning segment of $z''$.
Therefore $x\leq z$.
\bx

Note that for all $x,y\in G$:
\dy{d26b}{
$x\leq y$ if and only if $y^{-1}x\leq y^{-1}$.
}
Moreover, for all $x,y,z\in G$: 
\dy{d25a}{
if $x\leq y\leq z$ then $x^{-1}y\leq x^{-1}z, z^{-1}y\leq z^{-1}x$,
and $x\leq zy^{-1}x$.
}

This implies:

\prop{329a}{
For all $x,y,z\in G$, if $xy\leq z$ and $x|y$ then $x\leq zy^{-1}$.
}

\pf
Since $x\leq xy\leq z$, by \rf{d25a} we have
$x\leq z(xy)^{-1}x=zy^{-1}$.
\bx

In fact, the partial order $\leq$ yields a lattice if we add to $G$ an element $\infty$ at infinity.
First consider the following algorithm.
For any $x\in G$ let $\fff(x)$ denote the set of symbols $\alpha$ with $\alpha\leq x$.

For any two reduced words $x$ and $y$ the algorithm is as follows:
\dy{d7b}{
Grow a reduced word $z$ such that $zx'=x$ and $zy'=y$, where $zx'$ and $zy'$ are reduced words. 
Initially, $z:=\emptyset$.
If $z$ has been found, choose a symbol $\alpha\in\fff(x')\cap\fff(y')$, reset $z:=z\alpha$, remove the first occurrences of $\alpha$ from $x'$ and $y'$,  and iterate.
Stop if no such $a$ exists.
}
Clearly the final $z$ is reduced (as it is a beginning segment of a word arising by commuting $x$) and satisfies $z\leq x$ and $z\leq y$.
Note that $\fff(x')\cap\fff(y')=\emptyset$ if and only if the word $y'^{-1}x'$ is reduced.
Moreover:

\prop{pN5d}{
If $w\leq x$ and $w\leq y$ then $w\leq z$.
}

\pf
We may assume that $w$ is a reduced word.
Apply the algorithm \rf{d7b} to $w$ and $z$.
We end up with reduced words $vz'=z$ and $vw'=w$ such that $\fff(w')\cap\fff(z')=\emptyset$.
If $w'\equiv\emptyset$ then $w\leq z$, so assume $w'\not\equiv\emptyset$.
Let $x'$ and $y'$ be as found in \rf{d7b} applied to $x$ and $y$.
So $vz'x'=x$ and $vz'y'=y$ are reduced words.
Since $vw'=w\leq x=vz'x'$ the first symbol $a$ (say) of $w'$ belongs to $\fff(z'x')$. Similarly, $a$ belongs to $\fff(z'y')$.
Since $a\not\in\fff(z')$ it follows that $a$ and $z'$
are independent and that $a$ belongs to $\fff(x')\cap\fff(y')$.
This contradicts the construction of $z$.
\bx

It follows that $\leq$ forms a lattice on $G\cup\{\infty\}$, with $x\wedge y=z$ where $z$ is constructed as in \rf{d7b}.
Note that \rf{d7b} also gives an algorithm to test if $x\leq y$ for any two words.
Moreover:
\dy{d25c}{
for all $x,y\in G$: $x^{-1}(x\wedge y)\leq x^{-1}y$.
}
This follows from the fact that if $z,x'$ and $y'$ are as constructed in \rf{d7b} then $(x')^{-1}=x^{-1}(x\wedge y)$ and $(x')^{-1}y'=x^{-1}y$, while $(x')^{-1}y'$ is a reduced word.

One has $x\vee y$ is finite (i.e., belongs to $G$) if and only if there is  a $w\in G$ with $x\leq w$ and $y\leq w$.
The following proposition describes how to find $x\vee y$.
Let $x$ and $y$ be reduced.
Let $z,x'$ and $y'$ be as constructed in \rf{d7b}.

\prop{p5b}{
If $x'$ and $y'$ are independent then $x\vee y$ is finite and is equal to $zx'y'$.
Otherwise, $x\vee y=\infty$.
}

\pf
If $x'$ and $y'$ are independent, then $zx'y'=zy'x'$ is a reduced word and $x=zx'\leq zx'y'$ and $y=zy'\leq zy'x'$.

Now let $zx'\leq w$ and $zy'\leq w$ for some reduced word $w$.
Then $w$ can be commuted to $zx'w'$ and to $zy'w''$, where $x'w'$ and $y'w''$ can be commuted to each other.
Since $w''^{-1}y'^{-1}x'w'=1$ there exists a perfect matching $M$ satisfying \rf{d5f}.
As $x'w'$ and $y'w''$ are reduced, each symbol in the $w''^{-1}y'^{-1}$ part is matched with a symbol in the $x'w'$ part.
Since $\fff(x')\cap\fff(y')=\emptyset$, the symbols in the $y'^{-1}$ part are not matched with the symbols in the $x'$ part.
So $x'$ and $y'$ are independent and hence $zx'y'$ is a reduced word and $zx'y'\leq w$. 
\bx

This directly implies that for any $x,y\in G$ with $x\vee y$ finite:
\dy{dj6a}{
$x\vee y=x(x\wedge y)^{-1}y$ and $x\wedge y=x(x\vee y)^{-1}y$.
}
Note that:
\dy{d25b}{
if $x\leq z$ and $y\leq z$ then $x\vee y=z(z^{-1}x\wedge z^{-1}y)$ and $x\wedge y=z(z^{-1}x\vee z^{-1}y)$.
}
The reason is that by \rf{d25a} the function $w\mapsto z^{-1}w$ reverses the partial order on the set $\{w\in G \mid w\leq z\}$.
Hence $z^{-1}(x\vee y)=z^{-1}x\wedge z^{-1}y$ and $z^{-1}(x\wedge y)=z^{-1}x\vee z^{-1}y$ whenever $x,y\leq z$.

We also note the following:

\prop{p5c}{
Let $x_1,\ldots,x_t\in G$ be such that $x_i\vee x_j$ is finite for all $i,j$.
Then $x_1\vee\cdots\vee x_t$ is finite.
}

\pf
Let $x_1,\ldots,x_t$ be a counterexample with $t$ as small as possible.
So $x_1\vee\cdots\vee x_t=\infty$ and $t\geq 3$.
By the minimality of $t$ each pair from $x_1\vee x_4\vee\cdots\vee x_t, x_2\vee x_4\vee\cdots\vee x_t, x_3\vee x_4\vee\cdots\vee x_t$ has finite join.
If $t\geq 4$ this implies by the minimality of $t$ that $x_1\vee\ x_2\vee x_3\vee x_4\vee\cdots\vee x_t$ is finite, a contradiction.
So  $t=3$.

Assume we have chosen $x_1,x_2,x_3$ so that $|x_1|+|x_2|+|x_3|$ is as small as possible.
Then $x_1,x_2$ and $x_3$ are reduced nonempty words.
Write $x_1\equiv y_1\alpha, x_2\equiv y_2\beta, x_3\equiv y_3\gamma$ where $\alpha,\beta,\gamma$ are symbols.

By the minimality condition, $z:=y_1\vee y_2\vee y_3$ and $z\vee x_1$ are finite.
If $z\vee x_1=z$ then $x_1\vee x_2\vee x_3=y_1\vee x_2\vee x_3$ is finite by the minimality condition.
So $z\vee x_1\neq z$.
Since $y_1\leq z$ we know $z=y_1z'$ for some word $z'$ with $y_1|z'$.
Hence, since $x_1=y_1\alpha$ with $y_1|\alpha$, $z\vee x_1=za$.
Similarly $z\vee x_2=zb$ and $z\vee x_3=zc$.
Moreover $a\neq b$, since otherwise $x_1\vee x_2\vee x_3=za\vee zb\vee zc=za\vee zc=x_1\vee y_2\vee x_3$ is finite.

Since $za\vee zb=x_1\vee x_2\vee y_3$ is finite, $a$ and $b$ are independent.
Similarly  $a$ and $c$ are independent and  $b$ and $c$ are independent.
So $za,zb,zc\leq zabc=zbac=zcab$, and hence $x_1\vee x_2\vee x_3=za\vee zb\vee zc$ is finite.
\bx

The partial order $\leq$ is clearly not invariant under mappings $x\mapsto zx$ for $z\in G$.
The following formula expresses how $\wedge$ behaves under such an operation.

\prop{p31a}{
For all $x,y,z\in G$ one has $z^{-1}x\wedge z^{-1}y=z^{-1}((x\wedge y)\vee(x\wedge z)\vee(y\wedge z))$.
}

\pf
We may assume $x\wedge y\wedge z=1$ (we can replace $x$ by $(x\wedge y\wedge z)^{-1}x$, $y$ by $(x\wedge y\wedge z)^{-1}y$, and $z$ by $(x\wedge y\wedge z)^{-1}z$).
Let $u:=x\wedge y, v:= x\wedge z$ and $w:=y\wedge z$.
Since $u\wedge v=1$ and $u\vee v$ is finite, $u$ and $v$ are independent.
Similarly, $u$ and $w$ are independent and $v$ and $w$ are independent.
So $(x\wedge y)\vee(x\wedge z)\vee(y\wedge z)=uvw$.
Let $x':=(uv)^{-1}x, y':=(uw)^{-1}y$ and $z':=(vw)^{-1}z$.
So $uvx', uwy'$ and $vwz'$ represent $x,y$ and $z$ as reduced words.
Hence $(z')^{-1}w^{-1}ux'=(z')^{-1}uw^{-1}x'$ and $(z')^{-1}v^{-1}uy'=(z')^{-1}uv^{-1}y'$ represent $z^{-1}x$ and $z^{-1}y$ as reduced words.
Now $w^{-1}x'\wedge v^{-1}y'=1$ as $(w^{-1}x')^{-1}v^{-1}y'=(x')^{-1}v^{-1}wy'$ is a reduced word.
Hence $z^{-1}x\wedge z^{-1}y=(z')^{-1}u=z^{-1}uvw$.
\bx

It follows from Proposition \ref{p31a} that $(x\wedge y)\vee(x\wedge z)\vee(y\wedge z)$
is the unique element $w$ that is on the three shortest paths
(with respect to the distance function $\dist$) from $x$ to $y$,
$x$ to $z$, and $y$ to $z$.
So $(x\wedge y)\vee (x\wedge z)\vee(y\wedge z)$ is the `median'
in the sense of Sholander [15,\linebreak[0]17,\linebreak[0]16]
(cf. [2]).

\subsect{s12B}{Join-irreducible elements and partial distributivity of $G$}

It is helpful to see that each element of a free partially commutative
group has an underlying partial order --- extending the idea that the
symbols in a word in the free group are totally ordered.

Let $\alpha_1\cdots\alpha_t$ be a reduced word representing element $x$ of $G$.
Define a partial order $\preceq$ on $\{1,\ldots,t\}$ by:
\dy{423B}{
$i\preceq i'\Leftrightarrow$ there exist $i_0=i<i_1<\cdots<i_s=i'$ (with $s\geq 0$) such that $a_{i_{j-1}}$ and $a_{i_j}$ are not independent, for each $j=1,\ldots,s$,
}
for $i,i'\in\{1,\ldots,t\}$.

There is a one-to-one correspondence between linear extensions $i_1,\ldots,i_t$ of $1,\ldots,t$ (with respect to $\preceq$) and reduced words representing $x$.
Here we define a {\em linear extension} with respect to $\preceq$ as a permutation $i_1,\ldots,i_t$ of $1,\ldots,t$ such that if $i_j\preceq i_{j'}$ then $j\leq j'$, for all $j,j'\in\{1,\ldots,t\}$.

For any linear extension $i_1,\ldots,i_t$ the word $w:=\alpha_{i_1}\cdots\alpha_{i_t}$ is a reduced word representing $x$.
This follows from the fact that any linear extension can be obtained from $1,\ldots,t$ by iteratively choosing two consecutive elements $i_j,i_{j+1}$ with $i_j\not\preceq i_{j+1}$ and replacing them by $i_{j+1},i_j$.
So $w$ arises by commuting from $\alpha_1\cdots\alpha_t$.

Conversely, let $\alpha_{i_1}\cdots\alpha_{i_t}$ be a reduced word representing $x$, where $i_1,\ldots,i_t$ is a permutation of $1,\ldots,t$.
We may assume that we have chosen indices such that if $\alpha_{i_j}=\alpha_{i_{j'}}$ and $j<j'$ then $i_j<i_{j'}$.
Then $i_1,\ldots,i_t$ is a linear extension with respect to $\preceq$.
This follows iteratively from the fact that if $\alpha_{i_j}$ and $\alpha_{i_{j+1}}$ are independent, then $i_j\not\preceq i_{j+1}$; thus if $i_1,\ldots i_ji_{j+1}\cdots i_t$ is a linear extension, then so is $i_1\cdots i_{j+1}i_j\cdots i_t$.

Let $L_x$ denote the set of lower ideals of $(\{1,\ldots,t\},\preceq)$.
(A subset $I$ of $\{1,\ldots,t\}$ is a {\em lower ideal} if $i\in I$ and $j\preceq i$ implies $j\in I$.)
Then the partially ordered sets $(L_x,\subseteq)$ and $(\{y\in G\mid y\leq x\},\leq)$ are isomorphic.
The isomorphism is given as follows.
Let $y\leq x$.
Then there is a linear extension $i_1,\ldots,i_t$ and an $s\leq t$ such that $\alpha_{i_1}\cdots\alpha_{i_s}$ is a reduced word representing $y$. 
Then $\{i_1,\ldots,i_s\}$ is a lower ideal of $\{1,\ldots,t\}$.
(Indeed, if $i_j\preceq i_{j'}$ then $j\leq j'$; so if moreover $j'\leq s$ then $j\leq s$.)
Conversely, for each lower ideal $I$ of $\{1,\ldots,t\}$ there exists a linear extension $i_1,\ldots,i_t$ and an $s\leq t$ such that $I=\{i_1,\ldots,i_s\}$.
Then $\alpha_{i_1}\cdots\alpha_{i_s}$ is a reduced word representing an element $y\leq x$.
It is not difficult to see that this gives a one-to-one correspondence, bringing $\leq$ to $\subseteq$.

In particular it follows that (cf. [2])

\prop{423E}{
For each $x\in G$, the set $\{y\in G\mid y\leq x\}$, partially ordered by $\leq$, forms a distributive lattice.
}

\vspace{2mm}
(That is, if $a,b,c\leq x$ then $a\wedge(b\vee c)=(a\wedge b)\vee(a\wedge c)$ and $a\vee(b\wedge c)=(a\vee b)\wedge(a\vee c)$.)

\pf
This follows directly from the facts that the partially ordered sets $(\{y\in G\mid y\leq x\},\leq)$ and $(L_x,\subseteq)$ are isomorphic and that the collection $L_x$ is closed under taking unions and intersections.
\bx

(The whole lattice on $G\cup\{\infty\}$ is generally not distributive: if $a$ and $b$ are distinct generators then $a\wedge(b\vee b^{-1})=a\wedge\infty=a$ while $(a\wedge b)\vee(a\wedge b^{-1})=1\vee 1=1$.)

As a corollary we have:

\prop{423F}{
If $y\vee z$ is finite then $x\wedge(y\vee z)=(x\wedge y)\vee(x\wedge z)$.
If $x\vee y$ and $x\vee z$ are finite then $x\vee(y\wedge z)=(x\vee y)\wedge(x\vee z)$.
}

\pf
The first line follows from the fact that taking $x':=x\wedge(y\vee z)\leq y\vee z$ we have $x',y,z\leq y\vee z$, implying $x'=(x'\wedge y)\vee(x'\wedge z)=(x\wedge y)\vee(x\wedge z)$. 
The second line follows from $(x\vee y)\wedge(x\vee z)=((x\vee y)\wedge x)\vee((x\vee y)\wedge z)=x\vee((x\wedge z)\vee(y\wedge z))=x\vee(y\wedge z)$ (using the first line).
\bx

The correspondence also gives a correspondence for join-irreducible
elements of $G$.
An element $x$ of $G$ is called {\em join-irreducible} or a
{\em left-interval}
if $x\neq 1$ and if $y\vee z=x$ implies $y=x$ or $z=x$.
Clearly, each element $x$ of $G$ is equal to the join of all join-irreducible elements $y\leq x$.
Moreover, $x$ is join-irreducible if and only if the partial order $\preceq$ defined above has a unique maximum element.
In other words:

\prop{423A}{
$x$ is join-irreducible if and only if $x\neq 1$ and there exists a symbol $\alpha$ such that each reduced word representing $x$ has last symbol equal to $\alpha$.
}

\pf
To see necessity, let $x$ be join-irreducible.
Let $v$ and $w$ be reduced words representing $x$, with last symbols $\alpha$ and $\beta$, respectively.
If $\alpha\neq\beta$ then $x\alpha^{-1}\neq x\beta^{-1}$.
Hence $x\alpha^{-1}\vee x\beta^{-1}=x$, contradicting the fact that $x$ is join-irreducible.

To see sufficiency, let $x\neq 1$ and let there be a symbol $\alpha$ such that each reduced word representing $x$ has last symbol equal to $\alpha$.
Moreover, let $x=y\vee z$.
Then $x=x'y'z'$,
where $x':=y\wedge z, y':=x'^{-1}y$ and $z':=x'^{-1}z$, $x'|y'|z'$ and $y'$ and $z'$ independent.
So not both $y'$ and $z'$ can have $\alpha$ as last symbol, implying that $y'=1$ or $z'=1$.
Therefore $z=x$ or $y=x$.
\bx

Let $J_x$ denote the collection of join-irreducible elements $y\leq x$.
Then:

\prop{423G}{
The partially ordered sets $(J_x,\leq)$ and $(\{1,\ldots,t\},\preceq)$ are isomorphic.
}

\pf
By the above, the join-irreducible elements correspond to lower-ideals of $\{1,\ldots,t\}$ that have a unique maximum element.
So they are determined by their maximum element, and we have the required correspondence.
\bx

It follows that:

\prop{423H}{
For each $x\in G$, the number of join-irreducible elements $y\leq x$ is equal to $|x|$.
}

\pf
Directly from Proposition \ref{423G}.
\bx

\subsect{s12C}{Cyclically reduced words and periodicity}

Call an element $a$ of $G$ {\em cyclically reduced} if $a\wedge a^{-1}=1$.
So $a$ is cyclically reduced if and only if $a|a$.
Note that for each element $a$ of $G$ there exist unique $b,c$ such
that $a=bcb^{-1}$ with $b|c|b^{-1}$ and $c$ cyclically reduced.
In fact, $b=a\wedge a^{-1}$.

\prop{401a}{
Let $x,a\in G$, where $x$ is a left-interval with maximum symbol $\alpha$
and where $\alpha$ occurs in $a$.
Assume $x\leq ax$ and $b\wedge x=1$.
Write $a=bcb^{-1}$ with $b|c|b^{-1}$ and $c$ cyclically reduced.
Then $x$ and $b$ are independent and $x\leq c^t$ for some $t$.
}

\pf
Let $e:=a^{-1}\wedge x$.
Then $x\leq ax$ and $e|e^{-1}x$, and hence by Proposition \ref{329a}
$e\leq ae$. Since $ae\leq a$, we have
$e=e\wedge a=a^{-1}\wedge x\wedge a=b\wedge x=1$.
So $a\leq ax$.
Hence $b\leq a\leq ax$ and $x\leq ax$, and therefore $b\vee x<\infty$,
and so $b$ and $x$ are independent, and $\alpha$ occurs in $c$.

We next show that $x\leq c^t$ for some $t$ by induction on $|x|$.
We may assume that $x\not\leq c$.
Moreover $x\wedge c\neq 1$, since $x\leq cx$ and $c\leq cx$, and
hence if $x\wedge c=1$ then $x$ and $c$ are independent, but both contain
$\alpha$.

Let $x':=(x\wedge c)^{-1}x$ $c':=(x\wedge c)^{-1}c(x\wedge c)$
and $a':=(x\wedge c)^{-1}a(x\wedge c)=bc'b^{-1}$.
Then $x'\leq a'x'$, $\alpha$ occurs in $a'$ and $b\wedge x'=1$
(as $b$ and $x$ are independent).
So by induction we know that $x'\leq(c')^s$ for some $s$.
Hence $x\leq c^{s+1}$.
\bx

\prop{p3BB}{
Let $a$ be cyclically reduced and let $x\in G$.
Then $x\leq a^t$ for some $t$, if and only if
there exist $a_t\leq a_{t-1}\leq\cdots\leq a_2\leq a_1\leq a$ such that
$x=a_1a_2\cdots a_{t-1}a_t$ and $a_1|a_2|\cdots|a_{t-1}|a_t$ and
such that $a_i$ and $a^{-1}_{i-1}a$ are independent for each $i=2,\ldots,t$.}

\pf
Sufficiency being easy, we show necessity.
This is shown by induction on $t$, the case $t=1$ being trivial.
Let $a_1:=a\wedge x$ and $x':=a_1^{-1}x, b:=a_1^{-1}a$.
Since $a\leq a^t$ and $x\leq a^t$ we know that $a\vee x$ is finite, and hence $b$ and $x'$ are independent.
So $x'\leq a^{t-1}$ since $x'\leq ba^{t-1}$ and $b$ and $x'$ are independent.
So by induction there exist
$a_t\leq a_{t-1}\leq\cdots\leq a_2\leq a$ such that
$x'=a_2\cdots a_{t-1}a_t$ and $a_2|\cdots|a_{t-1}|a_t$
and such that $a_i$ and $a^{-1}_{i-1}a$ are independent for each $i=3,\ldots,t$.
Since $x'$ and $b$ are independent we know in fact that $a_2\leq a_1$
and that $a_2$ and $a_1^{-1}a$ are independent.
Thus we have the required $a_1,\ldots,a_t$.
\bx 

Let $a\in G$.
An element $d\in G$ is called a {\em component} of $a$ if
$d$ is a minimal element with the properties that $1\neq d\leq a$
and that $d$ and $d^{-1}a$ are independent.
So if $d_1,\ldots,d_n$ are the components of $a$ then
$a=d_1\cdots d_n$ where the $d_i$ are pairwise
independent.

\prop{401b}{
Let $a,x\in G$ be such that $x$ is a left-interval with maximum symbol $\alpha$,
such that $x\leq ax$ and
such that $x^{-1}ax$ contains $\alpha$.
Let $p$ be a natural number satisfying $p<\frac{|x|}{|a|}-|a|$.
Let $y$ be the component of $x^{-1}ax$ containing $\alpha$.
Then $x=ry^p$ where $r$ is a left-interval
with maximum symbol $\alpha$, and $r|y^p$.
}

\pf
I. Write $a=bcb^{-1}$ with $b|c|b^{-1}$ and $c$ cyclically reduced.
First assume that $b=1$.
Then by Proposition \ref{401a} $x\leq c^t$ for some $t$.
By Proposition \ref{p3BB}, there exist $c_t\leq c_{t-1}\leq\cdots\leq c_2\leq c_1\leq c$
such that $x=c_1c_2\cdots c_{t-1}c_t$ and $c_1|c_2|\cdots|c_{t-1}|c_t$
and such that $c_i$ and $c_{i-1}^{-1}c$ are independent for $i=2,\ldots,t$.
We may assume that $c_t\neq 1$.
As $t\geq |x|/|c|>|c|+p$,
there exists an $h\geq p+1$ such that $c_h=c_{h-1}$.
Hence $c_h$ and $c_h^{-1}c$ are independent.
Let $d:=c_h$ and $q:=c_h^{-1}c$.
Then $c_1=\cdots=c_h=d$.
Let $r:=c_{h+1}\cdots c_t$.
Then $x=d^h r$ and $r\leq x$ (by Proposition \ref{p3BB}).
Hence $x=r(x^{-1}dx)^h$ (since
$r^{-1}x=r^{-1}d^hr=(r^{-1}dr)^h=(x^{-1}dx)^h$.)
Moreover, $r|(x^{-1}dx)^h$.

Now $r$ is a left-interval with maximum symbol
$\alpha$, since $r$ is an end-segment of $x$.
Since $x$ and $q$ are independent, we have $x^{-1}ax=x^{-1}dxq$
with $x^{-1}dx$ and $q$ independent, and
therefore $x^{-1}dx=y$.
In particular, $x=ry^{p+1}$ with $r|y^{p+1}$.

II.
We now delete the assumption that $b=1$.
Let $x':=(b\wedge x)^{-1}x$, $b':=(b\wedge x)^{-1}b$ and
$a':=(b\wedge x)^{-1}a(b\wedge x)$.
By Proposition \ref{401a}, $b'$ and $x'$ are independent.
Moreover, $x'\leq cx'$, since
$x'\leq a'x'=b'c(b')^{-1}x'=b'cx'(b')^{-1}$.
Since moreover, $|x'|\geq|a'|^2+p|a'|+|a'|$,
we know by Section \ref{13me03a}
that $x'=r'(y')^{p+1}$, and $r'|(y')^{p+1}$,
where $y'$ is the component of $(x')^{-1}cx'$ containing $\alpha$,
and where $r$ is a left-interval with maximum symbol $\alpha$.

Let $y''$ be such that $(x')^{-1}cx'=y'y''$.
Since $b'$ and $x'$ are independent, also $b'$ and $y'$ are independent.
So $x^{-1}ax=(x')^{-1}a'x'=(x')^{-1}b'c(b')^{-1}x'=
b'(x')^{-1}cx'(b')^{-1}=b'(y'y'')(b')^{-1}=y'b'y''(b')^{-1})$
where $y'$ and $b'y''(b')^{-1}$ are independent.
So $y'=y$.

Now $b\wedge x|x'$ and $r'|y|y^p$
and hence $b\wedge x|r'|y|y^p$.
Now $r:=(b\wedge x)r'y$ is a left-interval with maximum symbol
$\alpha$.
For let $\beta\neq\alpha$ be a maximum symbol of $r$.
As $r'y$ is a left-interval with maximum symbol $\alpha$,
and as $r'|y$,
$\beta$ is a maximum symbol of $b\wedge x$ such that $\beta$
and $r'y$ are independent.
Then $\beta$ and $r'y^{p+1}$ are independent,
and hence $\beta$ is also a maximum symbol of $x$.
This contradicts the fact that $x$ is a left-interval with
maximum symbol $\alpha$.
\bx

\subsect{s26a}{Convex and closed sets}

We call a subset $H$ of $G$ {\em left-convex} if $H$ is nonempty and if $x,z\in H$ and $\dist(x,y)+\dist(y,z)=\dist(x,z)$ then $y\in H$.
Since the distance function is invariant under functions $x\mapsto yx$,
if $H$ is left-convex also $yH$ is left-convex for any $y\in G$.

\prop{p3a}{
A nonempty subset $H$ of $G$ is left-convex if and only if
\di{d3aa}{
\nr{i} if $x\leq y\leq z$ and $x,z\in H$ then $y\in H$;
\nrs{ii} if $x,y\in H$ then $x\wedge y\in H$ and, if $x\vee y$ is finite,
$x\vee y\in H$.
}
}

\pf
Necessity follows from the facts that if $x\leq y\leq z$ then $\dist(x,y)+\dist(y,z)=\dist(x,z)$, that $\dist(x,y)=\dist(x,x\wedge y)+\dist(x\wedge y,y)$ and that, if $x\vee y$ is finite then $\dist(x,y)=\dist(x,x\vee y)+\dist(x\vee y,y)$.

To see sufficiency, let $\dist(x,y)+\dist(y,z)=\dist(x,z)$ with $x,z\in H$.
We must show $y\in H$.
So $|x^{-1}y|+|y^{-1}z|=|x^{-1}z|$.
This implies $x^{-1}y|y^{-1}z$.
So $y^{-1}x\wedge y^{-1}z=1$.
Hence by Proposition \ref{p31a}, $(x\wedge y)\vee(x\wedge z)\vee(y\wedge z)=y$.
So $x\wedge z\leq x\wedge y\leq x$ and $x\wedge z\leq y\wedge z\leq z$.
Therefore, $x\wedge y, x\wedge z$ and $y\wedge z$ belong to $H$ and hence $y$ belongs to $H$.
\bx

In particular, each left-convex set has a unique minimum element.
We call a subset $H$ of $G$ {\em right-convex} if $H^{-1}$ is left-convex,
and {\em convex} if $H$ is both left- and right-convex.
(As usual, $H^{-1}:=\{x^{-1}\mid x\in H\}$.)
We call a subset $H$ of $G$ {\em left-closed} if $H$ is left-convex
and $1\in H$.
It is {\em right-closed} if $H^{-1}$ is left-closed, and {\em closed}
if $H$ is both left- and right-closed.

\prop{p3b}{
A nonempty subset $H$ of $G$ is left-closed if and only if 
\di{d3ba}{
\nr{i} if $y\leq x$ and $x\in H$ then $y\in H$;
\nrs{ii} if $x,y\in H$ and $x\vee y$ is finite then $x\vee y\in H$.
}
}

\pf
Directly from Proposition \ref{p3a}.
\bx

\prop{p27d}{
If $H$ is left-closed then for any $x\in G$: $x$ belongs to $H$ if and only if each left-interval of $x$ belongs to $H$.
If $H$ is closed then for any $x\in G$: $x$ belongs to $H$ if and only if each interval of $x$ belongs to $H$.
}

\pf
Let $H$ be left-closed.
If $x\in H$ and $y$ is a left-interval of $x$ then $y\leq x$ and hence $y\in H$.
The converse follows from the fact that $x=\bigvee\{y\mid y$ left-interval of $x\}$.

Let $H$ be closed.
If $x\in H$ and $y$ is an interval of $x$ then $y$ is a segment of $x$ and hence $y\in H$.
Conversely, suppose that each interval of $x$ belongs to $H$.
Then each right-interval $z$ of $x$ belongs to $H$ (since each left-interval $w$ of $z$ belongs to $H$ as $w$ is an interval of $x$).
So by the first statement of the proposition applied to $x^{-1}$ and $H^{-1}$ we have $x\in H$.
\bx

Clearly, the intersection of any number of left-convex sets is again left-convex.
Moreover, left-convex sets satisfy the following `Helly-type' property:

\prop{pj4a}{
Let $H_1,\ldots,H_t$ be left-convex sets with $H_i\cap H_j\neq\emptyset$ for all $i,j=1,\ldots,t$.
Then $H_1\cap\cdots\cap H_t\neq\emptyset$.
}

\pf
Suppose not. Choose a counterexample with $t$ minimal.
If $t\geq 4$ then each two of $H_1\cap H_4\cap\cdots\cap H_t, H_2\cap H_4\cap\cdots\cap H_t, H_3\cap H_4\cap\cdots\cap H_t$ have a nonempty intersection (as it is an intersection of $t-1$ sets from $H_1,\ldots,H_t$).
So these three sets have a nonempty intersection, a contradiction.

So $t=3$.
Choose $x\in H_1\cap H_2, y\in H_1\cap H_3, z\in H_2\cap H_3$.
Without loss of generality, $z=1$ (as we can replace $H_1,H_2,H_3$ by $z^{-1}H_1,z^{-1}H_2,z^{-1}H_3$).
Now $x\wedge y\in H_1\cap H_2\cap H_3$.
\bx

Proposition \ref{pj4a} implies the following.
As usual, define $HH':=\{xx'\mid x\in H,x'\in H'\}$.

\prop{p1a}{
Let $H$ be left-convex and let $H'$ and $H''$ be right-convex, with $H'\cap H''\neq\emptyset$.
Then $HH'\cap HH''=H(H'\cap H'')$.
}

\pf
Clearly $HH'\cap HH''\supseteq H(H'\cap H'')$.
To see the reverse inclusion, let $x\in HH'\cap HH''$.
So $x^{-1}H\cap H'^{-1}\neq\emptyset$ and $x^{-1}H\cap H''^{-1}\neq\emptyset$.
Since also $H'\cap H''\neq\emptyset$, by Proposition \ref{pj4a}, $x^{-1}H\cap (H'\cap H'')^{-1}\neq\emptyset$.
Hence $x\in H(H'\cap H'')$.
\bx

For any left-convex $H$ and $x\in G$, there is a unique element
$y$ in $H$ `closest' to $x$, that is, one minimizing $\dist(x,y)$.
To see this we may assume $H$ is left-closed.
Then by \rf{d3ba} $y$ is the largest element in $H$ satisfying $y\leq x$.

We denote this element by $\clos_H(x)$.
Then

\prop{p7g}{
Let $H$ be left-convex and let $H'$ be right-closed.
Let $x\in G$ and $y:=\clos_H(x)$.
Then $x\in HH'$ if and only if $y^{-1}x\in H'$.
}

\pf
Sufficiency being trivial, we prove necessity.
We may assume that $y=1$ (since replacing $H$ by $y^{-1}H$ and $x$ by $y^{-1}x$ does not modify the assertion).
Hence $H$ is left-closed.
Let $x\in HH'$ and choose $w\in H,w'\in H'$ such that $x=ww'$ and such that $|w|+|w'|$ is as small as possible.
We may assume that $w$ and $w'$ are reduced words.
Then the word $ww'$ is reduced again (since otherwise we could decrease $|w|+|w'|$).
So $w\leq x$ and hence $w=1$.
Therefore, $x=w'\in H'\subseteq HH'$.
\bx

Proposition \ref{p7g} implies that if $H$ is left-closed and $H'$ is right-closed, and if we can test in polynomial time whether any given word $x$ belongs to $H$ and to $H'$, then we can also test in polynomial time if any given word $y$ belongs to $HH'$.
We first find $z=\clos_H(x)$.
This can be done as follows: if we have a reduced word $x'x''=x$ with $x'\in H$, find a symbol $a\in\fff(x'')$ such that $x'a\in H$; reset $x':=x'a$, delete the first $a$ from $x''$, and iterate.
If no such $a$ exists we set $z:=x'$.
Now by Proposition \ref{p7g} $x\in HH'$ if and only if $y^{-1}x\in H'$. 

Similarly, if $H$ is left-convex and $H'$ is right-convex and if we can 
test in polynomial time if any word belongs to $H$ and to $H'$ and moreover
we know at least one word $w$ in $H$ and at least one word $w'$ in $H'$,
then we can test if any given word $x$  belongs to $HH'$: we just test
if $w^{-1}xw'^{-1}$ belongs to $(w^{-1}H)(H'w'^{-1})$.
Note that $w^{-1}H$ is left-closed and that $H'w'^{-1}$ is right-closed.

Proposition \ref{p7g} also implies:

\prop{p7f}{
Let $H$ be left-convex and let $H'$ be closed.
Then $HH'$ is left-convex.
}

\pf
We may assume that $H$ is left-closed.
We show that $HH'$ is left-closed.
Let $x\in HH'$ and $y\leq x$.
Let $u:=\clos_H(x)$. So by Proposition \ref{p7g}, $u^{-1}x\in H'$.
Moreover $(u\wedge y)^{-1}y\leq u^{-1}x$ (since $u\leq u\vee y\leq x$ and $u\vee y=u(u\wedge y)^{-1}y$).
This implies $(u\wedge y)^{-1}y\in H'$ and hence, as $u\wedge y\in H$, $y\in HH'$.

Next let $x,y\in HH'$ with $x\vee y$ finite.
Let $u:=\clos_H(x)$ and $v:=\clos_H(y)$.
So $u^{-1}x$ and $v^{-1}y$ belong to $H'$.
Now $(u\vee v)^{-1}(x\vee v)=(u\vee v)^{-1}(x\vee(u\vee v))=(x\wedge(u\vee v))^{-1}x=u^{-1}x\in H'$ (since $x\wedge(u\vee v)=u$ as $u=\clos_H(x)$).
Hence $x\vee v\in(u\vee v)H'$.
Similarly $y\vee u\in(u\vee v)H'$.
Hence $x\vee y=(x\vee v)\vee(y\vee u)\in(u\vee v)H'\subseteq HH'$.
\bx

For any $x\in G$ define
\dy{d27g}{
$H^{\uparrow}_x:=\{y\in G\mid y\geq x\}$ and $H^{\downarrow}_x:=\{y\in G\mid y\leq x\}$.
}
It is not difficult to see that $H^{\uparrow}_x$ and $H^{\downarrow}_x$ are left-convex.
 
\prop{p31e}{
Let $H$ be closed and let $x,y,z\in G$.
Then $x^{-1}yz$ belongs to $H$ if and only if: {\rm (i)} $\exists u\leq x~\exists w\leq z: u^{-1}yw\in H$; {\rm (ii)} $\exists u\leq x~\exists w\geq z: u^{-1}yw\in H$; {\rm (iii)} $\exists u\geq x~\exists w\leq z: u^{-1}yw\in H$; {\rm (iv)} $\exists u\geq x~\exists w\geq z: u^{-1}yw\in H$.
}

\pf
Necessity being trivial we show sufficiency.
Assertion (i) means $y\in H^{\downarrow}_xHH^{\downarrow-1}_z$, and assertion (ii) means $y\in H^{\downarrow}_xHH^{\uparrow-1}_z$.
Hence by Proposition \ref{p1a}, $y\in H^{\downarrow}_xH(H^{\downarrow-1}_z\cap H^{\uparrow-1}_z)= H^{\downarrow}_xHz^{-1}$.
Similarly, $y\in H^{\uparrow}_xHz^{-1}$.
Again by Proposition \ref{p1a}, $y\in (H^{\downarrow}_x\cap H^{\uparrow}_x)Hz^{-1}=xHz^{-1}$.
Therefore $x^{-1}yz\in H$.
\bx

\prop{p3CC}{
Let $H$ be closed and let $x,a$ be such that $x^{-1}ax\in H$.
Then there exists a $y$ such that $y^{-1}ay\in H$ and such that $y\leq a^t$ for some $t$.
}

\pf
By induction on $|x|$.
If $x\wedge ax\neq 1$, let $\alpha$ be a symbol satisfying $\alpha\leq x\wedge ax$.
Let $x':=\alpha^{-1}x$ and $a':=\alpha^{-1}a\alpha$.
Then $(x')^{-1}a'x'\in H$ and hence by induction there exists a $y'$ such that $(y')^{-1}a'y'\in H$ and such that $y'\leq (a')^s$ for some $s$.
Then $y:=\alpha y'$ satisfies $y^{-1}ay=(y')^{-1}a'y'\in H$ and $y\leq\alpha(a')^s\leq a^{s+1}$.

Similarly, if $x^{-1}\wedge (a^{-1}x^{-1})\neq 1$.
Hence we may assume that $x^{-1}|a|x$.
But then $a$ is a segment of $x^{-1}ax$, and hence $a\in H$.
This implies that we can take $y:=1$.
\bx

This implies:

\prop{p3DD}{
Let $H$ be closed and let $x,a$ be such that $x^{-1}a^sx\in H$ for some $s$.
Then there exists a $y$ such that $y^{-1}a^sy\in H$ and such that $y\leq a^{|a|}$.}

\pf
First assume that $a$ is cyclically reduced.
By Proposition \ref{p3CC} we may assume that $x\leq a^t$ for some $t$.
We choose $x$ and $t$ such that $t$ is minimal. 
We show that $t\leq |a|$.
Assume $t>|a|$.
By Proposition \ref{p3BB} there exist $1\leq a_t\leq a_{t-1}\leq\cdots\leq a_2\leq a_1\leq a$ such that $x=a_1a_2\cdots a_{t-1} a_t$ and
$a_1|a_2|\cdots|a_{t-1}|a_t$
and such that $a_i$ and $a_{i-1}^{-1}a$ are independent for each $i=2,\ldots,t$.
By the minimality of $t$, $a_t\neq 1$.
As $t>|a|$ there exist an $i\in\{2,\ldots,t\}$ such that $a_{i-1}=a_i$.
Then $a_i$ and $a_i^{-1}a$ are independent.
So $a_iaa_i^{-1}=a$.
Hence for $x':=a_i^{-1}x$ we have $(x')^{-1}a^sx'=x^{-1}a^sx\in H$ and $x'\leq a^{t-1}$, since $x'=a_1a_2\cdots a_{i-1}a_{i+1}\cdots a_t$.
This contradicts the minimality of $t$.

If $a$ is not cyclically reduced, let $a=bcb^{-1}$ with
$b|c|b^{-1}$ and $c$ cyclically reduced. 
Then for $x':=b^{-1}x$ we know that $(x')^{-1}c^sx'$ belongs to $H$,
and hence there exists a $y'$ such that $(y')^{-1}c^sy'\in H$ and such that $|y'|\leq c^{|c|}$.
Then for $y:=by'$ we have $y^{-1}a^sy\in H$ and $y\leq a^{|c|}\leq a^{|a|}$.
\bx

\sect{13me03b}{The cohomology feasibility problem}

\subsect{s3}{The cohomology feasibility problem}

Let $D=(V,A)$ be a directed graph and let $G$ be a group.
Two functions $\phi,\psi:A\to G$ are called
{\em cohomologous} if there exists a function $p:V\to G$ such that $\psi(a)=p(u)^{-1}\phi(a)p(w)$ for each arc $a=(u,w)$.
One directly checks that this gives an equivalence relation.

Consider the following {\em cohomology feasibility problem}:
\di{a3}{
\item[given:] a directed graph $D=(V,A)$, a group $G$, a function
$\phi:A\to G$, and for each $a\in A$, a 
subset $H(a)$ of $G$;
\items{find:} a function $\psi:A\to G$ such that
$\psi$ is cohomologous to $\phi$ and such that $\psi(a)\in H(a)$
for each $a\in A$.
}

We give a polynomial-time algorithm for this problem in case $G$ is an
free partially commutative group and each $H(a)$ is a closed set.

In the algorithm it is not required that the $H(a)$ are given explicitly.
It suffices to be given an algorithm that tests for any $a$ and any word
$x$ whether or not $x$ belongs to $H(a)$.
(So $H(a)$ might be infinite.)
The running time of the algorithm for the cohomology feasibility
problem is bounded by a polynomial
in $n:=|V|$, $\sigma:=\max\{|\phi(a)|\mid a\in A\}$, and $\tau$, where $\tau$ is the
maximum time needed to test membership 
of $x$ in $H(a)$ for any given arc $a$ and any given word $x$ of length bounded by a polynomial in $n$ and $\sigma$.
(The number $k$ of generators can be bounded by $n\sigma$, since we may assume that all generators occur among the $\phi(a)$.)

Note that, by the definition of cohomologous, equivalent to finding a $\psi$ as in \rf{a3}, is finding a
function $f:V\to G$ satisfying:
\dy{a3a}{
$f(u)^{-1}\phi(a)f(w)\in H(a)$ for each arc $a=(u,w)$.
}
We call such a function $f$ {\em feasible}.

Note that if $f$ is feasible and $P$ is an $s-t$ path, then $f(s)^{-1}\phi(P)f(t)\in H(P)$.
Here for any path $P=a_1\cdots a_m$ we use the following definitions:
\dy{a2}{
$\phi(P):=\phi(a_1)\cdots\phi(a_m)$, \\
$H(P):=H(a_1)\cdots H(a_m)$,
}
where $\phi(a^{-1}):=\phi(a)^{-1}$ and $H(a^{-1})=H(a)^{-1}$.

This gives an obvious necessary (but not sufficient) condition for
problem \rf{a3} having a solution:
\dy{14a}{
for each cycle $P$ there exists an $x\in G$ such that
$x^{-1}\phi(P)x$ belongs to $H(P)$.
}

\subsect{s3a}{Pre-feasible functions}

Let $D=(V,A)$ be a directed graph, let $G$ be a group, let $\phi:A\to G$ and
for each $a\in A$, let $H(a)$ be a closed subset of $G$.

We call a function $f:V\to G$ {\em pre-feasible} if for each arc $a=(u,w)$ of $D$ there exist $x\geq f(u)$ and $z\leq f(w)$ such that $x^{-1}\phi(a)z\in H(a)$.
Clearly, each feasible function is pre-feasible.
There is a trivial pre-feasible function $f$, defined by $f(v):=1$ for each $v\in V$.

The collection of pre-feasible functions is closed under certain operations on the set $G_M^V$ of all functions $f:V\to G_M$.
This set can be partially ordered by: $f\leq g$ if and only if $f(v)\leq g(v)$ for each $v\in V$.
Then $G_M^V$ forms a lattice if we add an element $\infty$ at infinity.
Let $\wedge$ and $\vee$ denote meet and join.

\prop{P25a}{
Let $f_1$ and $f_2$ be pre-feasible functions.
Then $f_1\wedge f_2$ and, if $f_1\vee f_2<\infty$, $f_1\vee f_2$ are pre-feasible again.
}

\pf
To see that $f_1\wedge f_2$ is pre-feasible, choose an arc $a=(u,w)$.
Since $f_1$ is pre-feasible, $\phi(a)\in H_{f_1(u)}^{\uparrow}H(a)(H_{f_1(w)}^{\downarrow})^{-1}\subseteq H_{f_1(u)\wedge f_2(u)}^{\uparrow}H(a)(H_{f_1(w)}^{\downarrow})^{-1}$.
Similarly, $\phi(a)\in H_{f_1(u)\wedge f_2(u)}^{\uparrow}H(a)(H_{f_2(w)}^{\downarrow})^{-1}$.
So by Proposition \ref{p1a}, 
\de{422a}{
\hspace{-16pt}\phi(a)\in H_{f_1(u)\wedge f_2(u)}^{\uparrow}H(a)(H_{f_1(w)}^{\downarrow}\cap H_{f_2(w)}^{\downarrow})^{-1}=H_{f_1(u)\wedge f_2(u)}^{\uparrow}H(a)(H_{f_1(w)\wedge f_2(w)}^{\downarrow})^{-1}.
}
This means that there exist $x\geq f_1(u)\wedge f_2(u)$ and $z\leq f_1(w)\wedge f_2(w)$ such that $x^{-1}\phi(a)z\in H(a)$.

The fact that $f_1\vee f_2$ is pre-feasible if it is not $\infty$, is shown similarly.
\bx

It follows that for each function $f:V\to G$ there is a unique smallest pre-feasible function $\bar{f}\geq f$, provided that there exists at least one pre-feasible function $g\geq f$.
If no such $g$ exists we set $\bar{f}:=\infty$.
Note that $\overline{f\vee g}=\bar{f}\vee\bar{g}$ for any two functions $f,g$ with
$f\vee g$ finite.

\subsect{s4}{A subroutine finding $\bar{f}$}

Let input $D=(V,A),\phi,H$ for the cohomology feasibility problem be given.
We describe a polynomial-time subroutine that outputs $\bar{f}$ for any given function $f$, under the assumption that \rf{14a} holds.

For any arc $a=(u,w)$ and any $x\in G$ let $\beta_a(x)$ be the smallest element $z$ in $G$ such that there exists an $x'\geq x$ with $(x')^{-1}\phi(a)z\in H(a)$.
This is unique, as $z$ is the minimum element in the left-convex set $\phi(a)^{-1}H^{\uparrow}_xH(a)$.
Note that for any $f:V\to G$ one has:
\dy{d14A}{
$f$ is pre-feasible if and only if $\beta_a(f(u))\leq f(w)$ for each arc $a=(u,w)$.
}

For any given $x$ we can determine $\beta_a(x)$ in polynomial time if we can test in
polynomial time if any given word belongs to $H(a)$ (as $\beta_a(x)$
is the minimal element of $\phi(a)^{-1}H^{\uparrow}_xH(a)$).

\vspace{4mm}
\noindent
{\bf Subroutine to find $\bar{f}$:} 
If $f$ is pre-feasible, output $\bar{f}:=f$.
Otherwise, choose an arc $a=(u,w)$ such that $\beta_a(f(u))\not\leq f(w)$.
If $f(w)\vee\beta_a(f(u))=\infty$, output $\bar{f}:=\infty$.
Otherwise reset $f(w):=f(w)\vee\beta_a(f(u))$, and start anew.

\vspace{2mm}
\prop{p8a}{
The output in the subroutine is correct.
}

\pf
Clearly, if $f(w)\vee\beta_a(f(u))=\infty$ then $\bar{f}=\infty$.
If $f(w)\vee\beta_a(f(u))<\infty$, let $f'$ denote the reset function.
Then $f\leq f'$.
Moreover, if $\bar{f}$ is finite, then $f\leq f'\leq\bar{f}$, since $f'(w)=f(w)\vee\beta_a(f(u))\leq\bar{f}(w)\vee\beta_a(\bar{f}(u))=\bar{f}(w)$, since $\bar{f}$ is pre-feasible.
\bx

\subsect{s5a}{Running time of the subroutine}

We show that after at most $2^{15}n^{9}k^{9}\sigma^{8}+2^2n^2k^2\rho$
iterations the subroutine gives an output, where:
\dy{d12B}{
$n:=|V|, \\
\sigma:=\max\{|\phi(a)|\mid a\in A\}, \\
\rho:=\max\{|f(v)|\mid v\in V\}.$
}

To this end we first make some observations and introduce some further
terminology.

\prop{p3c}{
$\beta_a(x\vee y)=\beta_a(x)\vee\beta_a(y)$ for all $x,y\in G$ with $x\vee y$ finite.
}

\pf
Clearly, for all $x,y$, if $x\leq y$ then $\beta_a(x)\leq\beta_a(y)$.
Equivalently, $\beta_a(x\vee y)\geq\beta_a(x)\vee\beta_a(y)$ for all $x,y$ with
$x\vee y$ finite.
To see the reverse inequality, let $u:=\beta_a(x)$ and $v:=\beta_a(y)$.
Since $x'^{-1}\phi(a)u\in H(a)$ for some $x'\geq x$, $H^{\uparrow}_x\cap\phi(a)H^{\downarrow}_{u\vee v}H(a)^{-1}\neq\emptyset$ (as $x'$ belongs to it).
Similarly, $H^{\uparrow}_y\cap\phi(a)H^{\downarrow}_{u\vee v}H(a)^{-1}\neq\emptyset$.
Since also $H^{\uparrow}_x\cap H^{\uparrow}_y=H^{\uparrow}_{x\vee y}\neq\emptyset$, by Proposition \ref{pj4a} $H^{\uparrow}_{x\vee y}\cap\phi(a)H^{\downarrow}_{u\vee v}H(a)^{-1}\neq\emptyset$.
Hence $\beta_a(x\vee y)\leq u\vee v$. 
\bx

Define for each path $P$ in $D$ and each $x\in G$, $\beta_P(x)$ inductively
by: $\beta_{\emptyset}(x):=x$ and $\beta_{Pa}(x):=\beta_a(\beta_P(x))$.
Inductively it follows from Proposition \ref{p3c} that $\beta_P(x\vee y)=\beta_P(x)\vee\beta_P(y)$ for all $x,y\in G$ with $x\vee y$ finite.
Moreover:

\prop{p9b}{
For each path $P$ in $D$, each $x\in G$ and each $y\in H(P)$ we have $\beta_P(x)\leq\phi(P)^{-1}xy$.
}

\pf
If $P=\emptyset$, the assertion is trivial.
If $P=a$, then for $z:=\phi(a)^{-1}xy$ one has $x^{-1}\phi(a)z=y\in H(a)$, and hence $\beta_a(x)\leq z=\phi(a)^{-1}xy$.

Consider next a path $Pa$ and let $y=y'y''\in H(Pa)$, with $y'\in H(P)$ and $y''\in H(a)$.
Then by induction, 
\dez{
\beta_{Pa}(x)=\beta_a(\beta_P(x))\leq\beta_a(\phi(P)^{-1}xy')
\leq\phi(a)^{-1}(\phi(P)^{-1}xy')y''=\phi(Pa)^{-1}xy.
}
\bx

We introduce the following further structure.
At each iteration $t$ of the subroutine we maintain a collection
$\Pi_t$ of paths.
Let $f_t$ denote the function $f$ as it is after $t$ iterations.
We first set $\Pi_0:=\{\emptyset\}$.
If at iteration $t$ we choose arc $a=(u,w)$ and put $f_t(w):=f_{t-1}(w)\vee\beta_a(f_{t-1}(u))$,
then for each left interval $x\leq\beta_a(f_{t-1}(u))$
satisfying $x\not\leq f_{t-1}(w)$ we choose a left-interval
$y\leq f_{t-1}(u)$ such that $x\leq\beta_a(y)$,
and we set $P_{w,x}:=P_{u,y}a$.
(Such a $y$ can be chosen by Proposition \ref{p3c}.)
We add each such path to $\Pi_{t-1}$, thus obtaining $\Pi_t$.

Note that the collection $\Pi_t$ has the following property:
\dy{318A}{
for each vertex $v$ and each left-interval $x\leq f_t(v)$ there is a
vertex $r$ and an $r-v$ path $P_{v,x}\in\Pi_t$ such that $x\leq\beta_{P_{v,x}}(f(r))$.
}

\prop{p22c}{
The subroutine takes at most
$t:=2^{15}n^{9}k^{9}\sigma^{8}+2^2n^2k^2\rho$ iterations. 
}

\pf
Suppose we have performed $t$ iterations.
We first show:
\dy{d10A}{
$\Pi_t$ contains a directed path $P$ with
$T:=2^{11}n^{7}k^{7}\sigma^{6}$ arcs.
}
First note that $|\Pi_t|\geq 2^{15}n^{9}k^{9}\sigma^{8}+2^2n^2k^2\rho$,
since at each iteration at least one path is added.
For each path $P_{v,x}$ in $\Pi_t$, with starting vertex $r$ say,
there exists a left-interval $y\leq f_0(r)$ such that
$x\leq\beta_{P_{v,x}}(y)$. 
Since there are  $n$ vertices and $2k$ symbols,
there exist vertices $r,v$ and symbols
$\alpha,\beta$ such that there are at least
$t/(2nk)^2=2^{13}n^{7}k^{7}\sigma^{8}+\rho$ join-irreducible
elements $x$ with the following properties:
\di{823a}{
\nr{i} $P_{v,x}$ belongs to $\Pi_t$ and runs from $r$ to $v$;
\nrs{ii} the maximum symbol of $x$ is equal to $\alpha$;
\nrs{iii} $x\leq\beta_{P_{v,x}}(y)$ for some left-interval
$y\leq f_0(r)$ with maximum symbol $\beta$.
}
Let $X$ denote the collection of such $x$.
Since each $x\in X$ has maximum symbol $\alpha$ and satisfies
$x\leq f_t(v)$, the elements in $X$ form a chain (i.e., are totally ordered by $\leq$).
Let $w$ be the maximum element in $X$.
Then
\de{d10B}{
|w|\geq |X|\geq 2^{13}n^{7}k^{7}\sigma^{8}+\rho.
}
Let $m$ be the number of arcs in $P_{v,w}$.
Let $z$ be the largest left-interval with maximum symbol $\beta$ 
and satisfying $z\leq f_0(r)$.
Then by Proposition \ref{p9b},
$w\leq\beta_{P_{v,w}}(z)\leq\phi(P_{v,w})^{-1}z$; so
$|w|\leq|\phi(P_{v,w})|+|z|\leq m\sigma+\rho$.
Hence with \rf{d10B},
$m\geq 2^{13}n^{7}k^{7}\sigma^{7}$.
Since each beginning segment of a path in $\Pi_t$ again belongs to
$\Pi_t$ we have \rf{d10A}.

Let $P$ traverse vertices $v_0,v_1,\ldots,v_T$ in this order.
By construction of $\Pi_t$ we can find left-intervals
$y_0,y_1,\ldots,y_T$ such that $y_0\leq f_0(v_0)$ and $y_i\leq\beta_a(y_{i-1})$
where $a=(v_{i-1},v_i)$ (for $i=1,\ldots,T$).
So $P_{v_i,y_i}$ is the subpath of $P$ consisting of the first $i$ arcs
of $P$.

For each vertex $v$ of $D$ and each symbol $\alpha$ let $I_{v,\alpha}$ denote the set of indices $i\in\{0,\ldots,T\}$ such that $v_i=v$ and such that $y_i$ has maximum symbol $\alpha$.
Then there exists a vertex $w$ of $D$ and a symbol $\beta$ such that $|I_{w,\beta}|\geq T/2nk$.
Let $L$ be the largest index in $I_{w,\beta}$.
Since for all $i,i'\in I_{w,\beta}$ and $i<i'$ we have $y_i<y_{i'}$ (since $y_i$ has maximum symbol $\beta$ and satisfies $y_i\leq f_t(w)$ for each $i\in I_{w,\beta}$ and since $y_{i'}\not\leq y_i$), we know
\dy{d11F}{
$|y_L|\geq |I_{w,\beta}|\geq T/2nk=2^{10}n^6k^6\sigma^5$.
}

Let $M:=2^4n^3k^3\sigma^2$ and $N:=2^3n^2k^2\sigma^2$.
Since $N=M/2nk$, there exists a vertex $u$ and a symbol $\alpha$ such that $I_{u,\alpha}$ contains at least $N+1$ indices $i$ satisfying $L-M\leq i\leq L$.
Choose $N+1$ such indices $i_0<i_1<\cdots<i_N$.
Define 
\de{d11H}{
x_0:=y_{i_0},x_1:=y_{i_1},\ldots,x_N:=y_{i_N}.
}
For $j=1,\ldots,N$ let $C_j$ be the $u-u$ path $v_{i_0},v_{i_0+1},\ldots,v_{i_j-1},v_{i_j}$.
We show that $C_N$ violates \rf{14a}.
Note that
\de{d10D}{
|\phi(C_j)|\leq (i_j-i_0)\sigma\leq M\sigma=2^5n^3k^3\sigma^3
}
for each $j=1,\ldots,N$.

Since $x_j\leq f_t(u)$ and since $\alpha$ is the maximum symbol of $x_j$ for each $j=0,\ldots,N$, we know that $x_0<x_1<\cdots <x_N$.
Moreover, 
\dez{
x_0<x_j\leq\beta_{C_j}(x_0)\leq\phi(C_j)^{-1}x_0
}
for each $j=1,\ldots,N$ (by Proposition \ref{p9b}).

By Proposition \ref{p9b}, $y_L\leq\beta_Q(x_0)\leq\phi(Q)^{-1}x_0$, where $Q$ denotes the path $v_{i_0},v_{i_0+1},\ldots, \\ v_{L-1},v_L$.
Hence, taking $p:= (M\sigma)^2$,
\de{r1}{
|x_0|\geq|y_L|-|\phi(Q)|\geq 4(M\sigma)^3-M\sigma\geq 3(M\sigma)^3
>|\phi(C_j)|^2+p|\phi(C_j)|.
}
Write $\phi(C_j)^{-1}=b_jc_jb_j^{-1}$ with
$b_j|c_j|b_j^{-1}$ and $c_j$ cyclically reduced.
Let $y_j$ be the component of $x_0^{-1}\phi(C_j)^{-1}x_0$
containing $\alpha$.
Then by Proposition \ref{401b}, 
$y_j$ is cyclically reduced and 
\dy{401c}{
$x_0=r_jy_j^p$ and $r_j|y_j^p$,
}
such that $r_j$ and $y_j$ are left-intervals with maximum symbol $\alpha$.

Let $y_j$ have $m_j$ symbols $\alpha$ and
let $x_0$ have $m$ symbols $\alpha$.
Write $x_0=z_mz_{m-1}\cdots z_1$, where each $z_i$
is a left-interval with maximum symbol $\alpha$
and where $z_m|z_{m-1}|\cdots|z_2|z_1$.
(Such a decomposition is unique.)
By \rf{401c} we know that, for each $j=1,\ldots,N$,
$m\geq pm_j$ and that
$z_i=z_{i'}$ if $i\equiv i'\pmod{m_j}$
and $i,i'\leq pm_j$.
Hence,
for $m:=\gcd\{m_1,\ldots,m_N\}$,
$z_i=z_{i'}$ if $i\equiv i'\pmod{m}$ and $i,i'\leq pm$.

Let $a:=z_nz_{n-1}\cdots z_1$ and $n_j:=m_j/m$ for
$m=1,\ldots,N$.
Then $y_j=a^{n_j}$ for each $j=1,\ldots,N$.

Write $x_0^{-1}\phi(C_j)^{-1}x_0=y_jy'_j$
for some $y'_j$ such that $y_j$ and $y'_j$ are independent.
Since $x_0^{-1}x_j\leq x_0^{-1}\phi(C_j)^{-1}x_0$ and since
$x_0^{-1}x_j$ is a left-interval with maximum symbol $\alpha$,
we know that $x_0^{-1}x_j\leq y_j$.
As $y_j^{n_{j'}}=y_{j'}^{n_j}$ for all $j,j'$, it follows that
$x_0^{-1}x_j$ and $y'_{j'}$ are  independent, for each $j'$.

Moreover,
\de{d10G}{
n_1<n_2<\cdots<n_N.
}
For suppose that $n_{j+1}\leq n_j$ for some $j=1,\ldots,N-1$.
Let $C$ be the closed path satisfying $C_{j+1}=C_jC$.
Then $x_{j+1}\leq\beta_C(x_j)\leq\phi(C)^{-1}x_j$ and hence
\dy{d10H}{
$x_0^{-1}x_{j+1}\leq x_0^{-1}\beta_C(x_j)\leq x_0^{-1}\phi(C)^{-1}x_j=x_0^{-1}\phi(C_{j+1})^{-1}\phi(C_j)x_j  \\
= (x_0^{-1}\phi(C_{j+1})^{-1}x_0)(x_0^{-1}\phi(C_j)x_0)(x_0^{-1}x_j) \\
= a^{n_{j+1}}a^{-n_j}(x_0^{-1}x_j)((y'_{j+1})^{-1}y'_j)$,
}
where $a^{n_{j+1}}a^{-n_j}(x_0^{-1}x_j)|((y'_{j+1})^{-1}y'_j)$.
This implies $x_{j+1}\leq (x_0a^{n_{j+1}-n_j})(x_0^{-1}x_j)$, and hence
$|x_{j+1}|\leq|x_0a^{n_{j+1}-n_j}|+|x_0^{-1}x_j|\leq|x_0|+|x_0^{-1}x_j|=|x_j|$.
(The inequality $|x_0a^{n_{j+1}-n_j}|\leq|x_0|$ follows from the fact that $x_0=fa^{n_j-n_{j+1}}$ for some $f$ satisfying $f|a^{n_j-n_{j+1}}$,
since $0\leq m_j-m_{j+1}\leq m_j\leq M\sigma$.)
This contradicts the fact that $x_{j+1}>x_j$, thus showing \rf{d10G}.

Now $|a|\leq M\sigma/N=2nk\sigma$, since $|a^{n_N}|\leq M\sigma$ and $n_N\geq N$ (by \rf{d10G}).
By \rf{14a}, there exists an $x\in G$ such that $x^{-1}\phi(C_N)x\in H(C_N)$.
Hence there exists a $y\in G$ such that $y^{-1}a^{n_N}y\in H(C_N)^{-1}$.
By Proposition \ref{p3DD} we may assume that $y\leq a^{|a|}$.
Hence $a^{n_N-|a|}\in H(C_N)^{-1}$.
Now by Proposition \ref{p9b},
$x_N\leq\beta_{C_N}(x_0)\leq\phi(C_N)^{-1}x_0a^{n_N-|a|}$,
and hence $x_0^{-1}x_N\leq x_0^{-1}\beta_{C_N}(x_0)\leq x_0^{-1}\phi(C_N)^{-1}x_0a^{|a|-n_N}=a^{|a|}y'_N$, with $a^{|a|}|y'_N$.
So $x_0^{-1}x_N\leq a^{|a|}$.
Therefore $|x_0^{-1}x_N|\leq |a|^2\leq (2nk\sigma)^2$.
Since $x_0< x_1<\cdots<x_N$ we know $|x_0^{-1}x_N|\geq N$.
Hence $2^3n^2k^2\sigma^2=N\leq (2nk\sigma)^2$, a contradiction.
\bx 

Combining the propositions, we obtain that the running time of the subroutine is bounded by a polynomial in the input size and in the time needed to check membership of $H(a)$. 
Then:

\thm{t22a}{
The running time of the subroutine is bounded by a polynomial in $n, k, \sigma$ and $\tau$, where
$\tau$ is the maximum time necessary to test if any word of size at
most $\rho+2^{24}n^{12}k^{12}\sigma^{11}$ belongs to any $H(a)$.
}

\pf 
Since initially $|f(v)|\leq\rho$ for each vertex $v$, we have after $t$ iterations $|f_t(v)|\leq\rho+\sigma t$ for each vertex $v$.
Since we do at most $2^{24}n^{12}k^{12}\sigma^{10}$ iterations, the result follows.
\bx

\subsect{s5}{A polynomial-time algorithm for the cohomology feasibility problem for free partially commutative groups}

We now describe the algorithm for the cohomology feasibility problem for free partially commutative groups.
Let $D=(V,A)$ be a directed graph, let $G$ be a free partially commutative group, let $\phi:A\to G$ and let $H(a)$ be a closed subset of $G$, for each $a\in A$.
We assume that with each arc $a=(u,w)$ also $a^{-1}=(w,u)$ is an arc, with $\phi(a^{-1})=\phi(a)^{-1}$ and $H(a^{-1})=H(a)^{-1}$.

Let $\UU$ be the collection of all functions $f:V\to G$ such that for each arc $a=(u,w)$ there exist $x\geq f(u)$ and $z\geq f(w)$ satisfying $x^{-1}\phi(a)z\in H(a)$.
For any given function $f$ one can check in polynomial time whether $f$ belongs to $\UU$.
Trivially, if $f\in\UU$ and $g\leq f$ then $g\in\UU$.
Moreover:

\prop{p14A}{
Let $f_1,\ldots,f_t$ be functions such that $f_i\vee f_j\in\UU$ for all $i,j$.
Then $f:=f_1\vee\cdots\vee f_t\in\UU$.
}

\pf
Choose an arc $a=(u,w)$.
We must show that for each arc $a=(u,w)$, $\phi(a)$ belongs to $H_{f(u)}^{\uparrow}H(a)(H_{f(w)}^{\uparrow})^{-1}$.
Since $f_i\vee f_j\in\UU$ for all $i,j$, we know
\de{d14D}{
\phi(a)\in H^{\uparrow}_{f_i(u)}H(a)(H^{\uparrow}_{f_j(w)})^{-1}}
for all $i,j$.
Hence by Proposition \ref{p1a}
\dy{d14E}{
$\displaystyle\phi(a)\in\bigcap_i\bigcap_j(H^{\uparrow}_{f_i(u)}H(a)(H^{\uparrow}_{f_j(w)})^{-1}) \\
=(\bigcap_iH^{\uparrow}_{f_i(u)})H(a)(\bigcap_jH^{\uparrow}_{f_j(w)})^{-1}=
H_{f(u)}^{\uparrow}H(a)(H_{f(w)}^{\uparrow})^{-1}$
}
Here $i$ and $j$ range over $1,\ldots,t$.
\bx

Let $X$ be the set of pairs $(u,x)$ where $u\in V$ and where $x$ is a left-interval such that there exists an arc $a=(u,w)$ with $x\leq\phi(a)$.
So $X$ has size polynomially bounded by $n$ and $\sigma$.
For any $(u,x)\in X$, let $f_{u,x}$ be the function defined by
\dt{d14B}{rcl}{
$f_{u,x}(u)$ & $:=$ & $x$, \\
$f_{u,x}(v)$ & $:=$ & $1$ for all $v\neq u$.
}

Let $E$ be the set of pairs $\{(u,x),(w,z)\}$ from $X$ such that there exists an arc $a=(u,w)$ such that
\dy{d12K}{
for all $x',z'\in G$, if $(x')^{-1}\phi(a)z'\in H(a)$ then $x\leq x'$ or $z\leq z'$.
}
Note that this holds if and only if $\phi(a)\not\in\Delta_xH(a)\Delta_z^{-1}$, where for any left-interval $y$, $\Delta_y$ is the left-closed set $\{y'\mid y'\not\geq y\}$.
So \rf{d12K} can be tested in polynomial time by Propositions \ref{p7g} and \ref{p7f}.

Let $E'$ be the collection of all pairs $\{(v,x),(v',x')\}$ from $X$ such that the function $\bar{f}_{v,x}\vee\bar{f}_{v',x'}=\infty$, or is finite and does not belong to $\UU$ (possibly $(v,x)=(v',x')$).

Choose a subset $Y$ of $X$ such that $e\cap Y\neq\emptyset$ for each
$e\in E$ and such that $e\not\subseteq Y$ for each pair $e\in E'$.
This is a special case of the 2-satisfiability problem, and hence can be solved in polynomial time.

\prop{p12A}{
If no such $Y$ exists, there is no feasible function.
}

\pf
Suppose $f$ is a feasible function.
Then $Y:=\{(v,x)\in X\mid v\in V, x\leq f(v)\}$ would have the required properties.
\bx

If we find $Y$, define $f$ by:
\de{d12N}{
f(v):=\bigvee\{\bar{f}_{v,x}\mid (v,x)\in Y\}.
}

\prop{p12C}{
 $f$ is a feasible function.}

\pf
Since $\bar{f}_{v,x}\vee\bar{f}_{v',x'}<\infty$ for each pair $\{(v,x),(v',x')\}\subseteq Y$, $f<\infty$.
Moreover, $f$ is the join of a finite number of pre-feasible functions, and hence $f$ is pre-feasible.
So by Proposition \ref{p31e} it suffices to show that for each arc $a=(u,w)$: 
\di{d14C}{
\nr{i} there exist $x\geq f(u)$ and $z\geq f(w)$ such that $x^{-1}\phi(a)z\in H(a)$;
\nrs{ii} there exist $x\leq f(u)$ and $z\leq f(w)$ such that $x^{-1}\phi(a)z\in H(a)$.
}

To see \rf{d14C}(i), note that it is equivalent to: $f\in\UU$.
As $\bar{f}_{v,x}\vee\bar{f}_{v',x'}\in\UU$ for all
$(v,x),(v',x')\in Y$, Proposition \ref{p14A} gives $f\in\UU$.

To see \rf{d14C}(ii), note that it is equivalent to:
\de{d12Q}{
\phi(a)\in H^{\downarrow}_{f(u)}H(a)(H^{\downarrow}_{f(w)})^{-1}.
}
Suppose \rf{d12Q} does not hold.
Let $b$ be the largest element in $H^{\downarrow}_{f(u)}H(a)$ satisfying
$b\leq\phi(a)$.
So by Proposition \ref{p7g},
$b^{-1}\phi(a)\not\in(H^{\downarrow}_{f(w)})^{-1}$; that is,
$\phi(a^{-1})b\not\leq f(w)$.
Hence there exists a left-interval $z$ of $\phi(a^{-1})b$ such that
$z\not\leq f(w)$.
So $\phi(a^{-1})b\not\in\Delta_z$ and hence by Proposition \ref{p7g},
$\phi(a)\not\in H^{\downarrow}_{f(u)}H\Delta_z^{-1}$.
Note that since $b\leq\phi(a)$ and $z\leq\phi(a^{-1})b$ we have
$z\leq\phi(a^{-1})$ and hence $(w,z)\in X$.

Let $c$ be the largest element in $\Delta_zH^{-1}$ such that
$c\leq\phi(a^{-1})$.
By Proposition \ref{p7g}, $\phi(a)c\not\in H^{\downarrow}_{f(u)}$;
that is $\phi(a)c\not\leq f(u)$.
Hence there exists a left-interval $x$ of $\phi(a)c$ such that
$x\not\leq f(u)$.
Again, since $c\leq\phi(a^{-1})$ and $x\leq\phi(a)c$ we have
$x\leq\phi(a)$ and hence $(u,x)\in X$.
So $\phi(a)c^{-1}\not\in\Delta_x$ and hence by Proposition \ref{p7g}, $\phi(a)\not\in\Delta_xH\Delta_z^{-1}$.
So $\{(u,x),(w,z)\}\in E$ and hence $Y$ contains at least one of $(u,x),(w,z)$.
So $x\leq f(u)$ or $z\leq f(w)$, a contradiction.
\bx

Thus we have proved:

\thm{t29a}{
The cohomology feasibility problem for free partially commutative groups is solvable in polynomial time.}
\bx

\subsect{s5b}{The 2-satisfiability problem}

In the algorithm we use a polynomial-time algorithm for the
2-satisfiability problem.
Conversely, the 2-satisfiability problem can be seen as a special case
of the cohomology feasibility problem for free groups.
To see this, first note that any instance of the 2-satisfiability problem can be
described as one of solving a system of inequalities in $\{0,1\}$ variables $x_1,\ldots,x_n$ of the form:
\dy{f2}{
$x_i+x_j\geq 1$ for each $\{i,j\}\in E$, \\
$x_i+x_j\leq 1$ for each $\{i,j\}\in E'$,
}
where $E$ and $E'$ are given collections of pairs and singletons from $\{1,\ldots,n\}$.
(So we allow $i=j$ in \rf{f2}, yielding $2x_i\geq 1$ or $2x_i\leq 1$.)

Let $G$ be the free group generated by the elements $g$ and $h$.
Make a directed graph with vertices $v_1,\ldots,v_{n}$ and with arcs:
\di{29f}{
\nr{i} $a=(v_i,v_j)$, with $\phi(a):=ghg^{-1}$, for each $\{i,j\}\in E$;
\nrs{ii} $a=(v_i,v_j)$, with $\phi(a):=h$, for each $\{i,j\}\in E'$. 
}
Moreover, set $H(a):=\{w\in G\mid |w|\leq 2\}$ for each arc $a$.

Now the cohomology feasibility problem in this case is equivalent to solving \rf{f2} in $\{0,1\}$ variables.
Indeed, if $x_1,\ldots,x_n$ is a solution of \rf{f2} then define $p(v_i):=g$ if $x_i=1$ and $p(v_i):=1$ if $x_i=0$. 
Then $p$ is a feasible function.
Conversely, if $p$ is a feasible function,
define $x_i:=1$ if $p(v_i)\neq 1$ and the first symbol of $p(v_i)$ is equal to $g$, and $x_i:=0$ otherwise.
Then $x_1,\ldots,x_n$ is a solution of \rf{f2}.

\subsect{s6}{A good characterization}

One may derive from the algorithm a `good' characterization of the
feasibility of the cohomology feasibility problem for free partially commutative groups,
i.e., one showing that the problem belongs to NP$\cap$co-NP.
We use the following well-known characterization for the feasibility
of \rf{f2}.
Assume that for all $h,i,j,k\in\{1,\ldots,n\}$:
\dy{f3}{
if $\{h,i\}\in E', \{i,j\}\in E, \{j,k\}\in E'$ then $\{h,k\}\in E'$.
}
(Extending iteratively $E'$ by any such pair $\{h,k\}$ does not change
the set of solutions of \rf{f2}.)

Then \rf{f2} has a $\{0,1\}$ solution if and only if 
\dy{f7}{
there is no $\{i,j\}\in E$ such that both $\{i\}$ and $\{j\}$ belong
to $E'$.
}
(If \rf{f3} does not hold, we could describe this condition in terms
of pairs of `alternating' cycles in $E\cup E'$.
If we would require moreover that \rf{f3} holds with $E$ and $E'$
interchanged, the condition will be that $E\cap E'$ does not contain
any singleton.)

We may adapt the subroutine in such a way that for each input $D=(V,A), \phi:A\to G, H(a)$ $(a\in A)$, and $f:V\to G$, we have as output:
\di{d15A}{
\nr{i} function $\bar{f}<\infty$, or
\nrs{ii} a cycle $C$ violating \rf{14a}, or
\nrs{iii} vertices $u,v,w$ of $D$, a directed $u-v$ path $P$ and a directed $w-v$ path $Q$ such that $\beta_P(f(u))\vee\beta_Q(f(w))=\infty$.
}  

\thm{t1}{
Let be given a directed graph $D=(V,A)$, a free partially commutative group $G$, a function $\phi:A\to G$, and for each arc $a$, a closed subset $H(a)$ of $G$.
Then there exists a function $\psi:A\to G$ such 
that $\psi$ is cohomologous to $\phi$ and $\psi(a)\in H(a)$ for each arc $a$, if and only if
\dy{18d}{
for each vertex $u$ and each two $u-u$ paths $P,Q$ there
exists an $x\in G$ such that $x^{-1}\cdot\phi(P)\cdot x\in H(P)$ and $x^{-1}\cdot\phi(Q)\cdot x\in H(Q)$. 
}}

\pf
{\bf Necessity.}
Let $f$ be a feasible function.
Then for $x:=f(u)$ we have
$x^{-1}\cdot\phi(P)\cdot x\in H(P)$ and $x^{-1}\cdot\phi(Q)\cdot x\in H(Q)$.

\vv{3mm}
\noindent
{\bf Sufficiency.}
Let \rf{18d} be satisfied, and assume that there is no feasible function $f$; that is, by Section \ref{s5}.
Note that \rf{18d} implies \rf{14a}.

Let $X, E$ and $E'$ be defined as in Section \ref{s5}.
We first show that $E$ and $E'$ satisfy \rf{f3}.
Let $\{(s,w),(t,x)\}\in E', \{(t,x),(u,y)\}\in E, \{(u,y),(v,z)\}\in E'$.
Assume $\{(s,w),(v,z)\}\not\in E'$; that is, $f:=\bar{f}_{s,w}\vee\bar{f}_{v,z}$ is finite and belongs to $\UU$.
 By definition of $E$, $a=(t,u)$ is an arc, and, since $f(t)^{-1}\phi(a)f(u)\in H(a)$, $x\leq f(t)$ or $y\leq f(u)$.
By symmetry we may assume $x\leq f(t)$.
This implies that $f_{t,x}\leq f$. 
Therefore, $f_{s,w}\vee f_{t,x}\leq f$, implying $\bar{f}_{s,w}\vee\bar{f}_{t,x}\leq f$.
So $\bar{f}_{s,w}\vee\bar{f}_{t,x}$ is finite and belongs to $\UU$.
This contradicts the fact that $\{(s,w),(t,x)\}\in E'$.

Since there is no feasible function, there is no subset $Y$ of $X$ such that $e\cap Y\neq\emptyset$ for each $e\in E$ and such that $e\not\subseteq Y$ for each $e\in E'$.
As \rf{f3} is satisfied it implies that there exists an arc $a=(u,w)$ and left-intervals $x,z$ such that $\{(u,x),(w,z)\}\in E$ and such that $\{(u,x)\},\{(w,z)\}\in E'$.
Then $x\leq \phi(a)$ and $z\leq\phi(a^{-1})$.
Since $\{u,x\}\in E'$, we know that $\bar{f}_{u,x}=\infty$ or is finite and does not belong to $\UU$.
It implies that 
\di{d15B}{
\nr{i}there exist a vertex $v$ and two $u-v$ paths $P,P'$ such that $\beta_P(x)\vee\beta_{P'}(x)=\infty$, or
\nrs{ii} there exist an arc $b=(v,v')$, a $u-v$ path $P$ and a $u-v'$ path $P'$ such that there do not exist $y\geq\beta_P(x)$ and $y'\geq\beta_{P'}(x')$ satisfying $y^{-1}\phi(a)y\in H(a)$.
}
Let $C$ be the $u-u$ cycle $P(P')^{-1}$ if (i) holds, and let $C$ be the $u-u$ cycle $Pb(P')^{-1}$ if (ii) holds.
Then 
\dy{d15C}{
there do not exist $c\geq x$ and $c'\geq x$ such that $c^{-1}\phi(C)c'\in H(C)$.
}
To see this, assume such $c,c'$ do exist.
Suppose first that \rf{d15B}(i) holds. 
Since $c^{-1}\phi(C)c'\in H(C)=H(P)H(P')^{-1}$, there exists an $y\in G$ such that $h:=c^{-1}\phi(P)y\in H(P)$ and $h':=(c')^{-1}\phi(P')y\in H(P')$.
Hence $\beta_P(x)\leq\beta_P(c)\leq\phi(P^{-1})ch=y$, and similarly $\beta_{P'}(x)\leq y$.
So $\beta_P(x)\vee\beta_{P'}(x)\leq y$, contradicting \rf{d15B}(i).

Suppose next that \rf{d15B}(ii) holds.
Since $c^{-1}\phi(C)c'\in H(C)=H(P)H(b)H(P')^{-1}$, there exist $y,y'$ such that $h:=c^{-1}\phi(P)y\in H(P), h':=(c')^{-1}\phi(P')y'\in H(P')$ and $y^{-1}\phi(b)y'\in H(b)$.
Hence $\beta_P(x)\leq\beta_P(c)\leq\phi(P^{-1})ch=y$, and similarly $\beta_{P'}(x)\leq y'$.
This contradicts \rf{d15B}(ii).

Similarly, there exists a $w-w$ cycle $D$ satisfying  
\dy{d15D}{
there do not exist $d\geq z$ and $d'\geq z$ such that $d^{-1}\phi(D)d'\in H(D)$.
}

By \rf{18d}, there exists a $c$ such that $c^{-1}\phi(C)c\in H(C)$ and $c^{-1}\phi(aDa^{-1})c\in H(aDa^{-1})$.
Hence there exist $d,d'$ such that $c^{-1}\phi(a)d\in H(a), d^{-1}\phi(D)d'\in H(D)$ and $(d')^{-1}\phi(a^{-1})c\in H(a^{-1})$.
By \rf{d15C}, $c\not\geq x$.
Since $c^{-1}\phi(a)d\in (a)$ and $\{(u,x),(w,z)\}\in E$ we know $d\geq z$.
Similarly, $d'\geq z$, contradicting \rf{d15D}.
\bx

\opm{o2}
Condition \rf{18d} cannot be  relaxed to requiring that for each cycle $P$ there exists an $x\in G$ such that $x^{-1}\cdot\phi(P)\cdot x$ belongs to $H(P)$.
To see this, let $G$ be the free group generated by $g$ and $h$.
Let $D$ be the directed graph with one vertex $v$ and two loops, $a$ and $b$, attached at $v$.
Define $\phi(a):=h,H(a):=\{1,h,g,g^{-1},g^{-1}h,hg\}$ and $\phi(b):=ghg^{-1},H(b):=\{1,h,g,g^{-1},hg^{-1},gh\}$.
If $x^{-1}\cdot\phi(a)\cdot x\in H(a)$ then the first symbol of $x$ is not equal to $g$.
If $x^{-1}\cdot\phi(b)\cdot x^{-1}\in H(b)$  then the first symbol of $x$ is equal to $g$.
So there is no $x$ such that both hold.

On the other hand, for each path $P$ there is an $x$ such that $x^{-1}\cdot\phi(P)\cdot x\in H(P)$.
Indeed, for each $k\in\oZ$, $\phi(ab^k)\in H(ab^k)$ and $\phi(b^ka)\in H(b^ka)$.
It follows that if $P$ starts or ends with $a$ or $a^{-1}$, then $\phi(P)\in H(P)$.
Moreover, for each $k\in\oZ$, $g^{-1}\cdot\phi(a^kb)\cdot g\in H(a^kb)$ and $g^{-1}\cdot\phi(ba^k)\cdot g\in H(ba^k)$.
So if $P$ starts and ends with $b$ or $b^{-1}$ then $g^{-1}\cdot\phi(P)\cdot g\in H(P)$.
\bx
 
The fact that Theorem \ref{t1} is a good characterization relies on the
facts that if the cohomology feasibility problem for free partially commutative groups has a solution, it has one
of small size, and that if paths $P,Q$ violating \rf{18d} would exist,
there are such paths of polynomial length.
(Both facts follow from the polynomial-time solvability of the 
subroutine.)
We can check in polynomial time whether or not for given $u-u$ paths
$P$ and $Q$ there exists an $x\in G$ such that $x\cdot\phi(P)\cdot x^{-1}$ belongs to $H(P)$ and $x\cdot\phi(Q)\cdot x^{-1}$ belongs to $\phi(Q)$.
(By the closedness of $H(P)$ and $H(Q)$ we have to consider
for $x$ only beginning segments of $\phi(P),\phi(P)^{-1},\phi(Q),\phi(Q)^{-1}$.
The number of such candidates for $x$ is polynomially bounded.)

\subsect{s7}{$R$-cohomologous functions}

In order to obtain results about paths instead of circuits, we extend the notion of cohomologous functions  to `$R$-cohomologous' functions.
Again, let $D=(V,A)$ be a directed graph, and let $(G,\cdot)$ be a group.
Moreover, let $R\subseteq V$.
Then two functions $\phi,\psi:A\to G$ are called $R${\em -cohomologous} if there exists a function $p:V\to G$ such that
\di{20a}{
\nr{i} $p(v)=\emptyset$ for all $v\in R$;
\nrs{ii} $\psi(a)=p(u)\cdot\phi(a)\cdot p(w)^{-1}$ for each arc $a=(u,w)$.
}
Again this defines an equivalence relation.

(One easily checks that if each component of $D$ contains at least one
vertex in $R$,
then $\phi$ and $\psi$ are equivalent, if and only if $\phi(P)=\psi(P)$ for each $r-s$ path $P$ with $r,s\in R$.
If $D$ is connected and $R=\{r\}$, there is a one-to-one correspondence
between $R$-cohomology classes and homomorphisms
$\Phi:\pi(D)\to G$, given by
$\Phi(\langle P\rangle):=\phi(P)$ for any $r-r$ path $P$.
Here $\pi(D)$ denotes the fundamental group of $D$ with base point $r$,
and $\langle P\rangle$ denotes the homotopy class containing path $P$.
Note that $\pi(D)$ itself is a free group.
We will not use these observations in the sequel.)

Consider the $R$-{\em cohomology feasibility problem\/}:
\di{20b}{
\item[given:] a directed graph $D=(V,A)$, a subset $R$ of $V$, a function
$\phi:A\to G$, and for each $a\in A$, a subset $H(a)$ of $G$;
\items{find:} a function $\psi:A\to G$ such that
$\psi$ is $R$-cohomologous to $\phi$ and such that $\psi(a)\in H(a)$ for each $a\in A$.
}
So equivalent is finding a function $p:V\to G$ such that $p(v)=0$ for all $v\in R$ and $p(u)\cdot\phi(a)\cdot p(w)^{-1}\in H(a)$ for each arc $a=(u,w)$.

If $G$ is a free partially commutative group $G$ and each $H(a)$ is closed, we can reduce problem \rf{20b} easily to the cohomology feasibility problem for free partially commutative groups.
We just add a loop at each vertex $v\in R$, add a new generator $g_0$ to the set of generators, and define $\phi(a):=g_0$ and $H(a):=\{\emptyset,g_0\}$ for each new arc (loop) $a$.
Let $\tilde{D}, \tilde{\phi}$ and $\tilde{C}$ denote the modified input.
One easily checks that the cohomology feasibility problem for $\tilde{D},\tilde{\phi},\tilde{C}$ is equivalent to the $R$-cohomology problem for $D,\phi,C$.
(Indeed, any feasible potential $p$ for $\tilde{D},\tilde{\phi},\tilde{C}$ should satisfy $p(v)=\emptyset$ for all $v\in R$.)

Thus we have:

\thm{t29b}{
The $R$-cohomology feasibility problem for free partially commutative groups is solvable in polynomial time.}
\bx

We can also derive from Theorem \ref{t1} a good characterization:

\thm{t2}{
Let be given a directed graph $D=(V,A)$, a subset $R$ of $V$, a free partially commutative group $G$, a function $\phi:A\to G$, and for each arc $a$, a closed subset $H(a)$ of $G$.
Then there exists a function $\psi:A\to G$ such
that $\psi$ is $R$-cohomologous to $\phi$ and $\psi(a)\in H(a)$ for each arc $a$,
if and only if
\di{20c}{
\nr{i} for each $r-s$ path $P$ with $r,s\in R$ one has $\phi(P)\in H(P)$;
\nrs{ii} for each vertex $s$ and each two $s-s$ paths $P,Q$ there
exists an $x\in G$ such that $x\cdot\phi(P)\cdot x^{-1}$ belongs
to $H(P)$ and $x\cdot\phi(Q)\cdot x^{-1}$ belongs to $H(Q)$. 
}
}

\pf
Necessity being trivial, we show sufficiency.
We extend $D,C,\phi$ to $\tilde{D},\tilde{C},\tilde{\phi}$ as above.
It suffices to show that \rf{20c} implies \rf{18d} (with respect to $\tilde{D},\tilde{C},\tilde{\phi}$).

Let $P$ and $Q$ be two $s-s$ paths in $\tilde{D}$, for some $s\in V$.
We must show that
\dy{20d}{
there exists an $x\in\tilde{G}$ such that $x\cdot\tilde{\phi}(P)\cdot x^{-1}\in\tilde{C}(P)$ and $x\cdot\tilde{\phi}(Q)\cdot x^{-1}\in\tilde{C}(Q)$.
}

\vspace{2mm}
I. If $P$ and $Q$ do not traverse any of the new loops attached at the points in $R$, then \rf{20d} directly follows from \rf{20c}(ii).

\vspace{2mm}
II. If both $P$ and $Q$ traverse some of the new loops, we can write $P=P_0a_1P_1\cdots a_mP_m$ and $Q=Q_0b_1Q_1\cdots b_nQ_n$, where $a_1,\ldots,a_m$ and $b_1,\ldots,b_n$ are new loops, and $P_0,\ldots,P_m$ and $Q_0,\ldots,Q_n$ are paths in the original graph $D$.

By \rf{20c}(i), $\phi(P_i)\in H(P_i)$ for $i=1,\ldots,m-1$ and $\phi(Q_i)\in H(Q_i)$ for $i=1,\ldots,n-1$.
Moreover, by the construction of $\tilde{D},\tilde{C},\tilde{\phi}$, one has that $\tilde{\phi}(a_i)\in\tilde{C}(a_i)$ for $i=1,\ldots,m$ and $\tilde{\phi}(b_i)\in\tilde{C}(b_i)$ for $i=1,\ldots,n$.

Consider the `surpluses' $\sigma(P_0^{-1}), \sigma(P_m), \sigma(Q_0^{-1}), \sigma(Q_n)$.
Let $x$ be one of largest size. 
We show that 
\dy{29g}{
$\phi(T)\cdot x^{-1}\in H(T)$ for each $T\in\{P_0^{-1},P_m,Q_0^{-1},Q_n\}$.
}
Without loss of generality, $x=\sigma(P_0^{-1})$.
Let $y:=\beta(P_0^{-1})$. So $\phi(P_0^{-1})=yx$.
If $x$ is an end segment of $\phi(T)$, then trivially $\sigma(T)$ is end segment of$x$ (as $x$ is at least as large as $x$), and hence $\phi(T)\cdot x^{-1}$ belongs to $H(T)$.
If $x$ is not an end segment of $\phi(T)$ then the last symbol of $\phi(T)\cdot x^{-1}$ is equal to the last symbol of $x^{-1}$.
So $\phi(T)\cdot\phi(P_0)=\phi(T)\cdot (x^{-1}y^{-1})=\phi(T)\cdot x^{-1})y^{-1}$ belongs to $H(T)\cdot H(P_0)$ (since $TP_0$ is an $r-r'$ path with $r,r'\in R$).
As by definition $y$ is the largest beginning segment of $\phi(P_0^{-1})$ that belongs to $H(P_0^{-1})$, it follows that $\phi(T)\cdot x^{-1}$ must belong to $H(T)$.
This proves \rf{29g}.

It implies that $x\cdot\tilde{\phi}(P)\cdot x^{-1}=(x\cdot \phi(P_0))\cdot\tilde{\phi}(a_1)\cdot\phi(P_1)\cdot\cdots\cdot\tilde{\phi}(a_m)\cdot(\phi(P_m)\cdot x^{-1})$ belongs to $H(P_0)\cdot H(a_0)\cdot H(P_1)\cdot\cdots\cdot H(a_m)\cdot H(P_m)=H(P)$.
Similarly for $Q$, thus proving \rf{20d}.

\vspace{2mm}

III. If only one of $P$ and $Q$ traverses some of the new loops, we may assume that $P$ does so.
Write $P=P_0a_1P_1\cdots a_mP_m$ such that $a_1,\ldots,a_m$ are new loops and $P_0,P_1,\ldots,P_m$ are paths in $D$.
As in part II one shows that there exists an $x\in G$ such that $x\cdot\phi(P_0)\in H(P_0)$ and $\phi(P_m)\cdot x^{-1}\in H(P_m)$.
We may assume that $x=\emptyset$, i.e., $\phi(P_0)\in H(P_0)$ and $\phi(P_m)\in H(P_m)$.
(We can reset $\phi(a):=x\cdot\phi(a)$ for each arc $a$ with tail $s$ and $\phi(a):=\phi(a)\cdot x^{-1}$ for each arc $a$ with head $s$.)

By \rf{20c}(ii), there exists a beginning segment $u$ of $\phi(Q)$ such
that $u^{-1}\cdot\phi(Q)\cdot u$ belongs to $H(Q)$.
If both $\phi(P_0^{-1})\cdot u\in H(P_0^{-1})$ and $\phi(P_m)\cdot u\in H(P_m)$, then as in part II, $u^{-1}\cdot\phi(P)\cdot u\in H(P)$, and we have \rf{20d}.
So we may assume that this is not the case.
Hence the largest beginning segment $y$ of $\phi(Q)$ such that both $\phi(P_0^{-1})\cdot y\in H(P_0^{-1})$ and $\phi(P_m)\cdot y\in H(P_m)$, satisfies $y\neq\phi(Q)$.
Similarly, we may assume that the largest beginning segment $z$ of $\phi(Q)^{-1}$ such that both $\phi(P_0^{-1})\cdot z\in H(P_0^{-1})$ and $\phi(P_m)\cdot z\in H(P_m)$, satisfies $z\neq\phi(Q)^{-1}$.

By definition of $y$ and $z$ we can choose $T,U\in\{P_0^{-1},P_m\}$ such that $y$ is the largest beginning segment of $\phi(Q)$ with $\phi(T)\cdot y\in H(T)$ and $z$ is the largest beginning segment of $\phi(Q)^{-1}$ with $\phi(U)\cdot z\in H(U)$.

We may assume that $z=\emptyset$.
(We can reset $\phi(a):=z\cdot\phi(a)$ for each arc $s$ with tail $a$ and $\phi(a):=\phi(a)\cdot z^{-1}$ for each arc with head $s$.)

First assume $y=\emptyset$. 
Then $\phi(TQU^{-1})=\phi(T)\phi(Q)\phi(U)^{-1}$. 
(Note that $\phi(T)\cdot\phi(Q)=\phi(T)\phi(Q)$ since $\phi(T)\in H(T)$ and $y=\emptyset$ is the largest beginning segment of $\phi(Q)$ with $\phi(T)\cdot y\in H(T)$.
Similarly, $\phi(Q)\cdot\phi(U)^{-1}=\phi(Q)\phi(U)^{-1}$.)
Since $\phi(TQU^{-1})$ belongs to $H(TQU^{-1})$ by \rf{20c}(i), it follows that $\phi(Q)$ belongs to $H(Q)$.
So we can take $x=\emptyset$ in \rf{20d}.  

Next assume $y\neq\emptyset$.
Then $\phi(Q)$ and $\phi(Q)^{-1}$ do not have any nonempty common beginning segment.
(Otherwise there is a nonempty common beginning segment $y'$ of $y$ and $\phi(Q)^{-1}$.
Then $\phi(U)\cdot y'\in H(U)$, contradicting the fact that $z=\emptyset$ is the largest beginning segment of $\phi(Q)^{-1}$ satisfying $\phi(U)\cdot z\in H(U)$.)

Now let $t:=|y|$ and let $y_0,\ldots,y_t$ be all beginning segments of $y$, with $|y_i|=i$ for $i=0,\ldots,t$.
Note that for each $i=0,\ldots,t$ one has 
\dy{21a}{
$\phi(P_m)\cdot y_i\in H(P_m)$.
}
This follows from the fact that $\phi(P_m)\cdot y_i$ is a beginning segment of at least one of $\phi(P_m)$ and $\phi(P_m)\cdot y$, where both belong to $H(P_m)$.
Similarly,
\dy{21c}{
$\phi(P_0^{-1})\cdot y_i\in H(P_0^{-1})$.
}
As in part II of this proof, \rf{21a} and \rf{21c} imply
\dy{22a}{
$y_i^{-1}\cdot\phi(P)\cdot y_i\in H(P)$ for each $i=0,\ldots,t$.
}

We show that for at least one $i\in\{0,\ldots,t\}$ one has 
\dy{21b}{
$y_i^{-1}\cdot\phi(Q)\cdot y_i$ belongs to $H(Q)$.
}
Combining this with \rf{22a} gives \rf{20d}.

Suppose \rf{21b} does not hold.
Let $z_i:=y_i^{-1}\cdot\phi(Q)y_{i-1}$ for $i=1,\ldots,t$.
Then
\de{29h}{
\phi(TQ^{t+1}U^{-1})=(\phi(T)\cdot y)z_tz_{t-1}\cdots z_2 z_1\phi(Q)\phi(U)^{-1}
} 
(i.e., no cancellations except at the $\cdot$; this follows from the facts that $z=\emptyset$, that $\phi(Q)$ and $\phi(Q)^{-1}$ have no nonempty common beginning segment, and that $y$ is the largest beginning segment of $\phi(Q)$ such that $\phi(T)\cdot y$ belongs to $H(T)$).

Then the assumption that \rf{21b} does not hold for any $i=0,\ldots,t$ implies that $\phi(TQ^{t+1}U^{-1})$ does not belong to $H(TQ^{t+1}U^{-1})$, contradicting \rf{20c}(i).
\bx

\sectz{Directed graphs on surfaces and homologous functions}

\subsect{s8}{Directed graphs on surfaces and homologous functions}

An embedding of a directed graph $D=(V,A)$ in a compact orientable surface $S$
(with each face being an open disk), can be described by a collection of cycles (`faces') $C_1,\ldots,C_f$ such that for each arc $a$ of $D$, each of $a$ and $a^{-1}$ occurs exactly once in $C_1,\ldots,C_f$.
For our purposes, such a cycle collection is enough to perform the algorithms below.
(We assume that each face is an open disk, which assumption does not restrict the generality of our results.)

We can think of the cycles $C_1,\ldots,C_f$ as giving the {\em clockwise} orientation of the faces.
In this interpretation, the face that traverses $a$ in forward direction is at the right-hand side of $a$, and the face that traverses $a$ in backward direction is at the left-hand side of $a$.

Note that by Euler's formula, $|V|+f=|A|+2-2h$, where $h$ is the number of handles of the surface.
Below when fixing a surface, we in fact just fix $h$.

For any directed graph$D=(V,A)$ embedded on a compact orientable surface, the {\em dual} graph $D^*=(\FF,A^*)$ has vertex set the collection $\FF$ of faces of $D$, while for any arc $a$ of $D$ there is an arc $a^*$ of $D^*$ with as tail the face of $D$ at the right-hand side of $a$ and as head the face of $D$ at the left hand side. 
We define for any function $\phi$ on $A$ the function $\phi^*$ on $A^*$ by $\phi^*(a^*):=\phi(a)$ for each $a\in A$. 

If a directed graph $D=(V,A)$ is embedded on a compact orientable surface $S$, we can dualize the concept of cohomologous functions to `homologous' functions.

Denote by $\FF$ the collection of faces of $D$.
Let $(G,\cdot)$ be a group.
We call two function $\phi,\psi:A\to G$ {\em homologous} if there exists a function $p:\FF\to G$ such that for each arc $a$ we have $p(F)\cdot\phi(a)\cdot p(F')^{-1}=\psi(a)$, where $F$ and $F'$ are the faces at the right-hand side and left-hand side of $a$, respectively.

The relation to cohomology is direct: $\phi$ and $\psi$ are homologous (in $D$), if and only if $\phi^*$ and $\psi^*$ are cohomologous (in $D^*$).

It follows that the {\em homology feasibility problem}:
\di{29i}{
\item[given:] a directed graph $D=(V,A)$ embedded on a compact orientable surface $S$, a function $\phi:A\to G$, and for each $a\in A$, a subset $H(a)$ of $G$;
\items{find:} a function $\psi:A\to G$ such that $\psi$ is homologous to $\phi$ and such that $\psi(a)\in H(a)$ for each $a\in A$,
}
is solvable in polynomial time if $G$ is a free partially commutative group and each $H(a)$ is closed.

\subsect{s6b}{Circulations and cycle decompositions}

Let $D=(V,A)$ be a directed graph embedded on a compact orientable surface $S$, and let $(G,\cdot)$ be a group.
We call a function $\phi:A\to G$ a {\em circulation} if for each vertex $v$ of $D$ we have 
\de{d3a}{
\phi(a_1)^{\varepsilon(v,a_1)}\cdot\ldots\cdot\phi(a_m)^{\varepsilon(v,a_m)}=1
}
where $a_1,\ldots,a_m$ are the arcs incident with $v$, in clockwise order, and $\varepsilon(v,a_i):=+1$ if $a_i$ enters $v$, and $:=-1$ if $a_i$ leaves $v$.
(If $a_i$ is a loop at $v$ we should be more careful.)

So $\phi$ is a circulation, if and only if for each cycle $\pi$ bounding a face of $D^*$ one has $\phi^*(\pi)=1$.
Note that in this last characterization it is not necessary to restrict oneself to clockwise cycles.
Consider e.g. three arcs, $a,b,c$ entering $v$, in clockwise order, with $\phi(a)\cdot\phi(b)\cdot\phi(c)=1$.
Then $\phi(c)\cdot\phi(b)\cdot\phi(a)$ is generally not equal to 1.
However, for $\pi:=a^*b^*c^*$, both $\phi^*(\pi)$ and $\phi^*(\pi^{-1})$ are equal to 1.

It is easy to check that if $\phi$ is a circulation and $\psi$ is homologous to $\phi$, then $\psi$ is again a circulation.

If $G$ is a free group, any circulation $\phi:A\to G$ can be decomposed as follows.
Replace any arc $a$ of $D$ by $t:=|\phi(a)|$ parallel arcs $a_1,\ldots,a_t$ (from right to left), yielding the graph $D_{\phi}=(V,A_{\phi})$.
Define $\phi'(a_i):=\xi_i$, where $\xi_i$ is the $i$th symbol in $\phi(a)$, for $i=1,\ldots,t$.

Consider now any vertex $v$.
Since \rf{d3a} holds we can find a perfect matching on the arcs of $D_{\phi}$ incident with $v$ (more precisely, on $\{1,\ldots,m\}$) in such a way that:
\di{23f}{
\nr{i} for any matched pair $\{a,b\}$ we have $\phi'(a)^{\varepsilon(v,a)}=\phi'(b)^{-\varepsilon(v,b)}$;
\nrs{ii} if $\{a,b\}$ and $\{c,d\}$ are matched pairs, thenthe  path $a^{-\varepsilon(v,a)}b^{\varepsilon(v,b)}$ does not cross the path $c^{-\varepsilon(v,c)}d^{\varepsilon(v,d)} $ at $v$.
}
Combining all matched pairs, at all vertices, we obtain a decomposition of $A_{\phi}$ into a collection $\CC$ of cycles, which we call a {\em cycle decomposition} of $\phi$.
The cycles do not have any fixed end point (formally speaking, we identify all cyclic permutations of the cycle).
No cycle in $\CC$ crosses itself or any of the other cycles in $\CC$.

Each cycle $C$ in $\CC$ has associated with it a symbol $\xi(C)$ from $g_1,g_1^{-1},g_2,g_2^{-1},\ldots$, such that for each arc $a$ of $D_{\phi}$, if $C$ traverses $a$ in forward direction then $\phi'(a)=\xi(C)$ and if $C$ traverses $a$ in backward direction then $\phi'(a)=\xi(C)^{-1}$.
The collection $\CC$ together with the function $\xi:\CC\to\{g_1,g_1^{-1},g_2,g_2^{-1},\ldots,\}$ uniquely determine $\phi$ (but generally not conversely).

We may consider the cycles in $\CC$ as cycles in $D$ (rather than in $D_{\phi}$) if, for each arc $a$ of $D$, we keep track of the order (from right to left) in which the cycles in $\CC$ traverse $a$.

\subsect{s6ba}{Disjoint circulations}

Let $D=(V,A)$ be a directed graph embedded on a compact orientable surface.
We call a circulation $\phi:A\to G_{\infty}$ {\em simple and directed} if any cycle decomposition of $\phi$ consists of pairwise vertex-disjoint simple directed cycles.
(A cycle is {\em simple} if no vertex is traversed more than once (except for the end vertices).
A cycle is {\em directed} if it does not contain $a^{-1}$ for any arc $a$.)

Consider the problem:
\di{d3b}{
\item[given:] a directed graph $D=(V,A)$ embedded on a compact orientable surface $S$ and a circulation $\phi:A\to G_{\infty}$;
\items{find:} a simple and directed circulation $\psi$ homologous to $\phi$.
}

In order to show that this problem is solvable in polynomial time, we define for each directed graph $D$ embedded on a compact orientable surface $S$, the `extended' dual graph $D^+=(\FF,A^+)$ as the graph obtained from $D^*$ by adding in each face of $D^*$ all chords.
(So generally $D^+$ is not embeddable in $S$.)
More precisely, for each nonempty path $\pi$ on the boundary of any face of $D^*$, $D^+$ has an arc $a_{\pi}$; if $\pi$ is an $F-F'$ path, $a_{\pi}$ runs from $F$ to $F'$.
(Since each arc $a^*$ is such a path, $D^+$ contains $D^*$ as a subgraph.)

For any $\phi:A\to G$, where $G$ is a group, define $\phi^+:A^+\to G$ by
$\phi^+(a_{\pi}):=\phi^*(\pi)$ for each $a_{\pi}\in A^+$.

\thm{t3a}{
Problem \rf{d3b} is solvable in polynomial time.}

\pf
Define
\dy{d3h}{
$H(a_{\pi}):=\{1,g_1,g_2,\ldots\}$ if $\pi:=a^*$ for any arc $a^*$ of $D^*$, \\
$H(a_{\pi}):=\{1,g_1,g_1^{-1},g_2,g_2^{-1},\ldots\}$ for all other arcs $a_{\pi}$ of $D^+$.
}
By Theorem \ref{t29a} we can find in polynomial time a function $\vartheta:A^+\to G_{\infty}$ such that $\vartheta$ is cohomologous to $\phi^+$ and such that $\vartheta(a_{\pi})\in H(a_{\pi})$ for each arc $a_{\pi}$ of $D^+$.
Defining $\psi(a):=\vartheta(a^*)$ for each $a\in A$ gives a solution of \rf{d3b}.

Moreover, if \rf{d3b} has a solution $\psi$, then such a function $\vartheta$ exists, viz.\ $\vartheta:=\psi^+$, as one directly checks.
\bx

\subsect{730}{The torus}

Theorem \ref{t1} implies a good characterization for the feasibility of \rf{d3b}, in terms of closed curves on $S$.
It is related to the one given in [11], where for any undirected graph $G$ embedded on a compact surface $S$ and any set of pairwise disjoint simple closed curves $C_1,\ldots,C_k$ on $S$, it was characterized when there exist pairwise disjoint simple circuits $C'_1,\ldots,C'_k$ in $G$ such that $C'_i$ is freely homotopic to $C_i$ for $i=1,\ldots,k$. 
(Freely homotopic means that there is no `base point'.)
However, in the present paper we consider the homology relation, which is coarser than homotopy, and so the two characterization do not seem to follow from each other.
 
However, if $S$ is the torus, the two concepts coincide.
This case has been dealt with by Seymour Seymour [14] (cf. Ding, Schrijver, and Seymour [3]).

Let $S=S^1\times S^1$ be the torus, where $S^1$ is a closed curve.
Let $S_1$ be the closed curve $S^1\times\{1\}$ on $S$ (fixing some orientation).
Let $C_1,\ldots,C_k$ be pairwise disjoint simple closed curves on $S$, each being freely homotopic to $S_1$ or to $S_1^{-1}$, choosing indices such that $C_1,\ldots,C_k$ occur cyclically around the torus (when going from the left-hand side of $S_1$ to the right-hand side).
We let the {\em sign} of $C_1,\ldots,C_k$ to be the vector $x\in\{+1,-1\}^k$ where $x_i=+1$ if $C_i$ is freely homotopic to $S_1$, and $x_i:=-1$ if $C_i$ is freely homotopic to $S_1^{-1}$.

For each closed curve $L$ on $S$ let the {\em winding number} $w(L)$ be equal to the number of times $L$ crosses $S_1$ from right to left, minus the number of times $L$ crosses $S_1$ from left to right.

Let $D=(V,A)$ be a directed graph embedded on the torus $S$, and let $L$ be a closed curve on $S$ with $w(L)\geq 0$.
We say that $L$ {\em fits} $x\in\{+1,-1\}^k$ if $L$ traverses points $p_1,\ldots,p_{kw(L)}$, in this order, such that for each $j=1,\ldots,kw(L)$:
\di{d3g}{
\item[either (i)] $p_j\in V$, 
\items{or (ii)} $p_j$ is on some arc $a$ of $D$ such that $D$ crosses $a$ from right to left if $x_{j}=+1$, and $D$ crosses $a$ from left to right if $x_{j}=-1$,
}
taking indices of $x_j$ modulo $k$.
We derive the following theorem of Seymour [14]:

\thm{t3d}{
Let $D=(V,A)$ be a directed graph embedded on the torus $S$ and let $x\in\{+1,-1\}^k$.
Then $D$ contains pairwise disjoint simple directed circuits each being freely homotopic to $S_1$ or $S_1^{-1}$, with sign $x$, if and only if 
\dy{d3i}{
each closed curve $L$ on $S$ fits some cyclic permutation of $x$.
}}

\pf
Necessity being trivial, we show sufficiency.
Let $k\geq 1$ and let \rf{d3i} be satisfied.
This easily implies that $D$ has at least one (undirected) circuit $C$ that is a freely homotopic to $S_1$.
Let $G$ be the free group generated by $g_1,\ldots,g_k$ and let $z:=g_1^{x_1}\cdots g_k^{x_k}$.
Define for each arc $a$ of $D$, $\phi(a):=z$ if $a$ that is traversed by $C$ in forward direction, $\phi(a):=z^{-1}$ if $a$ that is traversed by $C$ in backward direction, and $\phi(a):=1$ otherwise.
Let $H(a_{\pi})$ be as in \rf{d3h}.)
Each path $P$ in $D^+$ corresponds in a natural way to a curve on $S$ which we also denote by $P$.

We show that for each face $F$ and any two $F-F$ paths $P,Q$ in $D^+$ there exists an $x\in G$ such that 
\dy{d3j}{
$x\cdot\phi^+(P)\cdot x^{-1}\in H(P)$ and $x\cdot\phi^+(Q)\cdot x^{-1}\in H(Q)$.
}
We may assume that $w(P)\geq 0$ and $w(Q)\geq 0$.
Note that $\phi^+(P)=z^{w(P)}$ and $\phi^+(Q)=z^{w(Q)}$.
Assume that such an $x$ does not exist.
Define $z_i:=g_{i+1}^{x_{i+1}}\cdots g_k^{x_k}$, for $i=0,\ldots,k$.
By assumption, for each $i=1,\ldots,k$ there exists an $R_i\in\{P,Q\}$ such that $z_i\cdot\phi^+(R_i)\cdot z_i^{-1}\not\in H(R_i)$.
Let $R:=R_kR_{k-1}\cdots R_1R_0$. So $w(R)=w(R_k)+\cdots+w(R_0)$ and $\phi^+(R)=z^{w(R)}$.
By \rf{d3i}, some cyclic permutation of $\phi^+(R)$ belongs to $H(R)$.
Hence $z^{w(R)-1}g_1^{x_1}$ belongs to $H(R)$.
Since $z^{w(R)-1}=\prod_{i=k}^1(z_iz^{w(R_i)}\cdot z_{i-1}^{-1})$ and $H(R)=\prod_{i=k}^1H(R_i)$, there exists an $i=1,\ldots,k$ such that $z_iz^{w(R_i)}\cdot z_i^{-1}$ belongs to $H(R_i)$, a contradiction.
This show that there exists an $x\in G$ satisfying \rf{d3j}.

Now by Theorem \ref{t3} there exists a function $\vartheta:A^+\to G$ cohomologous to $\phi^+$ such that $\vartheta(a_{\pi})\in H(a_{\pi})$ for each arc $a_{\pi}$ of $D^+$.
Define $\psi(a):=\vartheta(a_{a^*})$ for each arc $a$ of $D$.
Then $\psi$ is a simple and directed circulation homologous to $\phi$.
Since for each cycle $P$ in $D^*$ one has $\psi^*(P)=\phi*(P)$ it follows that any cycle decomposition of $\psi$ consists of directed circuits of the required type.
\bx

A stronger version, also given by Seymour [14], in which we prescribe for each arc $a$ of $D$ which of the directed circuits $C_1,\ldots,C_k$ are permitted to traverse $a$, can also be derived.
(To this end we restrict the $H(a_{a^*})$.)

\subsect{s11a}{$\RR$-homologous function, $\delta$-joins, and path decompositions}

Dual to $\RR$-cohomologous functions are $\RR$-homologous functions.
Let $D=(V,A)$ be a directed graph embedded on a compact orientable surface $S$, with face collection $\FF$, and let $\RR\subseteq \FF$.
Let $(G,\cdot)$ be a group.
We call two functions $\phi,\psi:A\to G$ {\em $\RR$-homologous} if there exists a function $p:\FF\to G$ such that $p(F)=1$ for each $F\in\RR$ and such that for each arc $a$ we have $p(F)\cdot\phi(a)\cdot p(F')^{-1}=\psi(a)$, where $F$ and $F'$ are the faces at the right-hand side and left-hand side of $a$, respectively.
Again it follows that the {\em $\RR$-homology feasibility problem for free partially commutative groups} is solvable in polynomial time.

An extension of the notion of circulation is the `$\delta$-join'.
Let $\delta:V\to G$ (a `demand function') be such that
\dy{d31b}{
each vertex $v$ with $\delta(v)\neq 1$ has degree one.
}
(That is, $v$ is incident with exactly one arc.)
Call a function $\phi:A\to G$ a $\delta${\em -join} if for each vertex $v$: 
\de{23e}{
\phi(a_1)^{\varepsilon(v,a_1)}\cdot\ldots\cdot\phi(a_m)^{\varepsilon(v,a_m)}=\delta(v),
}
where again $a_1,\ldots,a_m$ are the arcs incident with $v$ in clockwise order, and for any arc $a$ incident with $v$, $\varepsilon(v,a):=+1$ if $a$ leaves $v$ and $\varepsilon(v,a):=-1$ if $a$ enters $v$.
So if $\delta(v)=1$ for each vertex $v$, any $\delta$-join is a circulation.

Let $W$ be the set of vertices $v$ satisfying $\delta(v)\neq 1$, and let $\RR$ be the collection of faces incident with at least one vertex in $W$.
One directly checks:
\dy{29j}{
if $\phi$ is a $\delta$-join and $\psi$ is $\RR$-homologous to $\phi$, then $\psi$ is a $\delta$-join again.
}

An extension of the idea of the cycle decomposition of a circulation is that of a `path decomposition' of a $\delta$-join.
Let $G$ be a free group.
Again we make the directed graph $D_{\phi}=(V,A_{\phi})$.
At each vertex $v\not\in W$ we can find a matching as for circulations.
Combining all matched pairs we obtain a decomposition of $A_{\phi}$ into a collection $\PP$ of paths and cycles, which we call a {\em path decomposition} of $\phi$.
Each path has tail and head in $W$ (possibly the same vertex).
The cycles do not have any given fixed end point (formally speaking, we identify all cyclic permutations of the cycle).
All vertices traversed by the paths except for their ends, and all vertices traversed by the cycles, belong to $V\setminus W$.
Moreover, none of the paths and cycles crosses itself or any of the other paths and cycles.

Each path and cycle $P$ in $\PP$ has associated with it a symbol $\xi(P)$ from $\{g_1,g_1^{-1},g_2,g_2^{-1},\ldots\}$, such that for each arc $a$ of $D_{\phi}$, if $P$ traverses $a$ in forward direction then $\phi'(a)=\xi(P)$ and if $P$ traverses $a$ in backward direction then $\phi'(a)=\xi(P)^{-1}$.  
The pair $\PP,\xi$ determines $\phi$.

\subsect{s8a}{Enumerating homology types}

With the methods developed before we can find in polynomial time a $\delta$-join of given homology type.
In order to be able to consider all homology types of a certain restricted
size, we describe an enumeration. 
(A related enumeration was given in [12].)
   
We call a $\delta$-join $\phi$ {\em elementary} if $\phi$ has a path decomposition with paths only (all starting and ending in $W$).
We first consider the following problem for any $p$ and compact orientable
surface $S$:
\di{d31c}{
\item[given:] 
a directed graph $D=(V,A)$ embedded on $S$, with exactly $p$ 
faces, a natural number $m$, and a function $\delta:V\to G_{\infty}$ such that each vertex in the set $W:=\{v\mid \delta(v)\neq 1\}$ has degree 1;
\items{find:} all elementary $\delta$-joins $\phi$ with $|\phi(a)|\leq m$ for each arc $a$ not incident with any vertex in $W$.
}

\thm{t31a}{
For each fixed $p$ and compact orientable surface $S$, problem \rf{d31c}
is solvable in polynomial time.
}

\vn{2mm} [As input size we take $|V|+|A|+m+\sum_{v\in V}|\delta(v)|$.]

\pf
We may assume that $|V\setminus W|=1$.
To see this, consider any arc $a$ connecting two different vertices in $V\setminus W$.
Let $D'=(V',A')$ arise from $D$ by contracting $a$.
Let $\delta'(v):=\emptyset$ for the contracted vertex $v$, and let $\delta'$ coincide with $\delta$ on all other vertices.
Then for each $\delta'$-join $\phi'$ there is a unique $\delta$-join
$\phi$ such that $\phi|A'=\phi'$.
So any enumeration of $\delta'$-joins gives directly an enumeration of
$\delta$-joins.

Let $V=W\cup\{u\}$ for some vertex $u$.
Hence $D$ consists of one vertex $u$, with a number of oriented loops at $u$, and a number of arcs connecting $u$ with the vertices in $W$, each of degree one.
We may assume that each of the nonloops has tail in $W$ and head $u$.
Let $L$ denote the set of loops of $D$.
By Euler's formula, $|L|=p+2h-1$.
So when considering the arcs incident with $u$ in clockwise order, there are 
$2p+4h-2$ (possibly empty) consecutive groups of nonloops, separated
by loops.
Let the $j$th group, $W_j$ say, consist of the arcs $(w_{j,1},u),\ldots,(w_{j,t_j},u)$, in clockwise order.
Let $\tilde{D}$ arise from $D$ by identifying for each $j$ all vertices
$w_{j,1},\ldots,w_{j,t_j}$ to one vertex $w_j$ and identifying all parallel arcs $(w_j,u)$ arising.

Consider any elementary $\delta$-join with $|\phi(a)|\leq m$ for each
loop $a$ of $D$.
Let $P_1,\ldots,P_M$ form a path decomposition of $\phi$.
Define the {\em type} of a path $P_i$ as the path in $\tilde{D}$ obtained from $P$ by traversing the arcs in $\tilde{D}$ that are parallel to those in $P$.

Let $Q_1,\ldots,Q_K$ be all types of $P_1,\ldots,P_M$ that are different from $aa^{-1}$ for any arc $a$. 
We identify type $Q$ with $Q^{-1}$.
Now $\phi$ is completely determined by the $Q_i$, together with 
a word $y_i$ in $G$ associated to $Q_i$ (for each $i$) that forms the concatenation of the
symbols asociated with the $P_j$ of type $Q_i$ (in the appropriate order).

We show that we can choose the $Q_j$ with the associated words in a polynomially bounded number of ways (fixing $p$ and $h$).
First we show that $K$ is at most $9p+18h$.
Since $P_1,\ldots,P_M$ are pairwise noncrossing, we can decouple
the paths $Q_1,\ldots,Q_K$ at $u$, so as to obtain pairwise disjoint
paths $Q'_1,\ldots,Q'_K$ (disjoint except for their end points).
Since any two $Q_i$ are different, the graph $H$ with vertex $u$ and arcs (loops) $Q'_1,\ldots,Q'_K$ has no faces bounded by one or two edges, except for
the $p$ original faces of $D$.
So $3(f-p)\leq 2K$, where $f$ denotes the number of faces of $H$.
As $\tilde{D}$ has at most $(2p+4h-2)+1$ vertices, by Euler's formula we have, $K\leq 2p+4h-1+f-2+2h\leq 3p+6h+\frac{2}{3}K$.
Therefore $K\leq 9p+18h$.

For any $a,b\in L\cup L^{-1}$, let $m_{ab}$ be the number of times $ab$ or $b^{-1}a^{-1}$ occurs in $Q_1,\ldots,Q_K$ (counting multiplicities).
Since the $P_i$ are pairwise noncrossing, we can reconstruct $\{Q_1,\ldots,Q_K\}$ from the $m_{ab}$ (up to reversing a path).

Now $m_{ab}\leq m$ for all $a,b$, since $|\phi(a)|\leq m$ for each $a\in A$.
So in enumerating, we can choose for each pair $a,b\in L\cup L^{-1}$ a nonnegative integer $m_{ab}\leq m$.
(Since $|L|=p+2h-1$, there are at most $(m+1)^{(p+2h-1)^2}$ choices.)
For each choice, we try to construct $Q_1,\ldots,Q_K$ from the $m_{ab}$.
If we fail or if $K>9p+18h$, we go on to the next choice of the $m_{ab}$.
If we succeed and $K\leq 9p+18h$, we proceed as follows.

For each $i=1,\ldots,K$, if the first arc of $Q_i$ equals $(w_j,u)$, the word associated with $Q_i$ should be equal to $\bar{x}$, where $x$ is some segment of the word 
\de{24f}{
\delta(w_{j,1})\cdots\delta(w_{j,t_j}).
}
Here $\bar{x}$ denotes the word obtained from $x$ by cancelling iteratively all occurrences of $\xi\xi^{-1}$ and $\xi^{-1}\xi$. 
(In \rf{24f} we did not cancel occurrences of $\xi\xi^{-1}$ or $\xi^{-1}\xi$.
So word \rf{24f} need not be in $G_{\infty}$.)

Since $K\leq 9p+18h$ there are at most $(\sum_{i=1}^{n}|\delta(w_i)|)^{18p+36h}$ such segments, and we can consider all choices in polynomial time.
Combining all choices gives us a function on the arcs of $D$, that
is either an elementary $\delta$-join, or not.
This way we obtain all elementary $\delta$-joins.
\bx

Consider next the following problem for any compact orientable surface $S$:
\di{d29k}{
\item[given:] a directed graph $D=(V,A)$ embedded on $S$, a natural number $m$, and a function $\delta:V\to G_{\infty}$, such that each vertex $v$ in the set $W:=\{v\mid \delta(v)\neq\emptyset\}$ has degree one;
\items{find:} $\delta$-joins $\phi_1,\ldots,\phi_N$ such that each elementary $\delta$-join $\phi$ with $|\phi(a)|\leq m$ for each arc not incident with $W$, is $\RR$-homologous to at least one of $\phi_1,\ldots,\phi_N$, where$\RR$ is the collection of faces incident with at least one vertex in $W$.
}

\thm{t3}{
For each fixed $p$ and compact orientable surface $S$, problem \rf{d29k} is solvable in polynomial time when $|\RR|=p$.
}

\vspace{2mm}
\noindent
[Again we take as input size $|V|+|A|+m+\sum_{v\in V}|\delta(v)|$.]

\pf
Let $A'$ denote the set of arcs not incident with any vertex in $W$.
Delete iteratively arcs from $D$ that are incident with at least one
face not in $\RR$.
We end up with a graph $\tilde{D}=(V,\tilde{A})$ that has $p$ faces, and such that
each elementary $\delta$-join in $D$ with $|\phi(a)|\leq m$
for each arc $a\in A'$, is $\RR$-homologous
to some elementary $\delta$-join in $\tilde{D}$ with $|\phi(a)|\leq 2m|A|$
for each arc $a\in A'$.
Thus Theorem \ref{t31a} implies the required enumeration.
\bx

\subsect{s9}{Applications to disjoint paths and trees problems}

We apply the techniques described above to a number of disjoint paths and disjoint trees problems.

We first consider the following problem, for any fixed compact surface $S$ and any fixed $p$:
\di{24g}{
\item[{\rm given:}] a directed graph $D=(V,A)$ embedded on $S$ and pairs $(r_1,s_1),\ldots,(r_k,s_k)$ of vertices of $D$, with the property that there exist $p$ faces such that each of $r_1,s_1,\ldots,r_k,s_k$ is incident with at least one of these faces;
\items{find:} pairwise vertex-disjoint paths $P_1,\ldots,P_k$, where $P_i$ is an $r_i-s_i$ path ($i=1,\ldots,k$).
}

\thm{t4}{
For each fixed compact orientable surface $S$ and each fixed $p$, problem \rf{24g} is solvable in polynomial time.
}

\pf
We may assume that $r_1,s_1,\ldots,r_k,s_k$ all are distinct and have degree one.
Let $\RR$ be the collection of faces incident with at least one of $r_1,s_1,\ldots,r_k,s_k$.
Define $\delta(r_i):=g_i$ and $\delta(s_i):=g_i^{-1}$ for $i=1,\ldots,k$.
Moreover, define $\delta(v):=\emptyset$ for all other vertices $v$.

By Theorem \ref{t3} we can find in polynomial time (fixing $S$ and $p$) a list of $\delta$-joins $\phi_1,\ldots,\phi_N$ in $D$ such that each elementary $\delta$-join $\phi$ with $|\phi(a)|\leq 1$ for each arc $a$ not incident with $r_1,s_1,\ldots,r_k,s_k$, is $\RR$-homologous to at least one of the $\phi_i$.

Consider the extended dual graph $D^+$ of $D$ (cf. Section \ref{s6ba}).
Define for each arc $a_{\pi}$ of $D^+$:
\dy{29n}{
$H(a_{\pi}):=\{\emptyset,g_1,\ldots,g_k\}$ if $\pi=a^*$ for some $a\in A$; \\
$H(a_{\pi}):=\{\emptyset,g_1,g_1^{-1},\ldots,g_k,g_k^{-1}\}$ for all other $a_{\pi}$ .
}
By Theorem \rf{t29a} we can find in polynomial time a function $\vartheta$ that is $\RR$-cohomologous to $\phi_i^+$ in $D^+$, with $\vartheta(b)\in H(b)$ for each arc $b$ of $D^+$, provided that such a $\vartheta$ exists.
If we find one, define $\psi(a):=\vartheta(a^*)$, for each arc $a$ of $D$.
Then $\psi$ is a $\delta$-join in $D$ (as it is $\RR$-homologous to $\phi_i$), and any path decomposition of $\psi$ into paths of cycles contains pairwise disjoint paths $P_1,\ldots,P_k$ as required.

If for none of $i=1,\ldots,N$ we find such a $\vartheta$ we may conclude that problem \rf{24g} has no solution.
For suppose $P_1,\ldots,P_k$ is a solution.
Define $\phi(a):=g_i$ if $P_i$ traverses $a$ ($i=1,\ldots,k$) and $\phi(a):=1$ if $a$ is not traversed by any $P_1,\ldots,P_k$.
Since $\phi$ is an elementary $\delta$-join with $|\phi(a)|\leq 1$ for each arc $a$, there exists an $i\in\{1,\ldots,N\}$ such that $\phi$ and $\phi_i$ are $\RR$-homologous.
However, for this $i$, there exists a $\vartheta$ as above, viz.\ $\vartheta:=\phi^+$.
This contradicts our assumption.
\bx

A special case applies to \rf{d1}:

\cor{30a}{
For each fixed $k$, the $k$ disjoint paths problem for directed planar graphs is solvable in polynomial time.
}

\pf
Directly from Theorem \ref{t4}.
\bx

An extension of Theorem \ref{t4} applies to the following problem: 
\di{24k}{
\item[{\rm given:}] a directed graph $D=(V,A)$ embedded on $S$ and pairs $(r_1,S_1),\ldots,(r_k,S_k)$ with $r_1,\ldots,r_k\in V$ and $S_1,\ldots,S_k\subseteq V$, with the property that there exist $p$ faces such that each vertex in $\{r_1,\ldots,r_k\}\cup S_1\cup\cdots\cup S_k$ is incident with at least one of the faces;
\items{find:} pairwise vertex-disjoint rooted trees $T_1,\ldots,T_k$, where $T_i$ has root $r_i$ and covers $S_i$ ($i=1,\ldots,k$).
}

\thm{t5}{
For each fixed compact orientable surface $S$ and each fixed $p$, problem \rf{24k} is solvable in polynomial time.
}

\pf
We may assume that all vertices in $W:=\{r_1,\ldots,r_k\}\cup S_1\cup\cdots\cup S_k$ are distinct and have degree one.
Let $\RR$ be the collection of faces incident with at least one vertex in $W$.
Define $\delta(r_i):=g_i^{|S_i|}$ and $\delta(s):=g_i^{-1}$ for $s\in S_i$, for $i=1,\ldots,k$.
Moreover, define $\delta(v):=1$ for all other vertices $v$.

By Theorem \ref{t3} we can find in polynomial time (fixing $S$ and $p$) a list of $\delta$-joins $\phi_1,\ldots,\phi_N$ in $D$ such that each elementary $\delta$-join $\phi$ with $|\phi(a)|\leq |V|$ for each arc $a$ not incident with $r_1,s_1,\ldots,r_k,s_k$, is $\RR$-homologous to at least one of the $\phi_i$.

Again consider the extended dual graph $D^+$ of $D$.
Define for each arc $a_{\pi}$ of $D^+$:
\dy{29p}{
$H(a_{\pi}):=\{g_i^n\mid i=1,\ldots,k;n\in\oZ, n\geq 0\}$ if $\pi=a^*$ for some $a\in A$; \\
$H(a_{\pi}):=\{g_i^n\mid i=1,\ldots,k;n\in\oZ\}$ for all other $a_{\pi}$.
}
By Theorem \rf{t29a} we can find in polynomial time a function $\vartheta$ that is $\RR$-cohomologous to $\phi_i^+$ in $D^+$, with $\vartheta(b)\in H(b)$ for each arc $b$ of $D^+$, provided that such a $\vartheta$ exists.
If we find one, define $\psi(a):=\vartheta(a^*)$, for each arc $a$ of $D$.
Then $\psi$ is a $\delta$-join in $D$ (as it is $\RR$-homologous to $\phi_i$), and any path decomposition of $\psi$ into paths and cycles contains pairwise disjoint rooted trees $T_1,\ldots,T_k$ as required.

If for none of $i=1,\ldots,N$ we find such a $\vartheta$ we may conclude that problem \rf{24k} has no solution.
For suppose $T_1,\ldots,T_k$ is a solution.
Define $\phi(a):=g_i^l$ if $T_i$ contains $a$ and for $l$ vertices $s$ in $S_i$ the simple $r-s_i$ path in $T_i$ traverses $a$, and $\phi(a):=1$ if $a$ is not contained in any $T_1,\ldots,T_k$.
Since $\phi$ is an elementary $\delta$-join with $|\phi(a)|\leq |V|$ for each arc $a$, there exists an $i\in\{1,\ldots,N\}$ such that $\phi$ and $\phi_i$ are $\RR$-homologous.
However, for this $i$, there exists a $\vartheta$ as above, viz.\ $\vartheta:=\phi^+$.
This contradicts our assumption.
\bx
 
A further extension is to the following problem:

\di{24m}{
\item[{\rm given:}] a directed graph $D=(V,A)$ embedded on $S$, subsets $A_1,\ldots,A_k$ of $A$, and pairs $(r_1,S_1),\ldots,(r_k,S_k)$ with $r_1,\ldots,r_k\in V$ and $S_1,\ldots,S_k\subseteq V$, with the property that there exist $p$ faces such that each vertex in $\{r_1,\ldots,r_k\}\cup S_1\cup\cdots\cup S_k$ is incident with at least one of these faces;
\items{find:} pairwise vertex-disjoint rooted trees $T_1,\ldots,T_k$, where $T_i$ has root $r_i$, covers $S_i$ and contains arcs only in $A_i$ ($i=1,\ldots,k$).
}

\thm{t6}{
For each fixed compact orientable surface $S$ and each fixed $p$, problem \rf{24m} is solvable in polynomial time.
}

\pf
As before, now replacing the first line in \rf{29p} by: $H(a_{\pi}):=\{g_i^n\mid i=1,\ldots,k, a\in A_i; n\in\oZ,n\geq 0\}$ if $\pi=a^*$ for $a\in A$.
\bx

We do not see if our methods extend to compact {\em non}orientable surfaces.

\subsect{s10}{Other groups and the arc-disjoint case}

Our algorithms are based on the polynomial-time solvability of the cohomology feasibility problem for free groups.
It might be interesting to investigate in how far the method can be extended to other groups.
Especially, for which groups $(G,\cdot)$ and subsets $C\subseteq G$ is the following problem solvable in polynomial time:
\di{24p}{
\item[{\rm given:}] a directed graph $D=(V,A)$ and a function $\phi:A\to G$;
\items{find:} a function $\psi:A\to C$ cohomologous to $\phi$.
}

This might apply to the arc-disjoint case as follows.
It is unknown if the following problem is solvable in polynomial time or NP-complete for $k=2$: 
\di{24n}{
\item[{\rm given:}] a directed planar graph and vertices $r_1,s_1,\ldots,r_k,s_k$, 
\items{find:} find pairwise arc-disjoint paths $P_1,\ldots,P_k$, where $P_i$ is an $r_i-s_i$ path ($i=1,\ldots,k$).
}
If we do not require planarity the problem is NP-complete for $k=2$, as follows from the result of Fortune, Hopcroft, and Wyllie mentioned in section \ref{s1}.
(The vertex-disjoint case can be reduced to the arc-disjoint case.)
The complexity status of \rf{24n} is also unknown for the special case $k=2,r_1=s_2,s_1=r_2$.

Now if problem \rf{24p} is polynomial-time solvable for the group $G:=\oZ^2$, taking $C=\{(0,0),(1,0),$$(0,1)\}$, then problem \rf{24n} is solvable in polynomial time for $k=2$.
This can be seen with a method similar to the one described in the previous sections.

More generally, if the cohomology feasibility problem for free groups is solvable in polynomial time for the group $\oZ^k$, taking for $H(a)$ the set of all unit basis vectors together with the origin, then problem \rf{24n} is polynomial-time solvable for this $k$.
Note that $\oZ^k$ can be considered as the `free abelian group'; it is generated by $g_1,\ldots,g_k$, with relations $g_i\cdot g_j=g_j\cdot g_i$ for all $i,j=1,\ldots,k$.

This also implies that if the group itself is part of the input of problem \rf{24p} (given, e.g., by generators and relations), then the problem will be NP-hard. 
This follows from the NP-completeness of problem \rf{24p} for nonfixed $k$.
Note that the algorithm we described for the cohomology feasibility problem for free groups is polynomial-time also if we do not fix the number of generators.

In a sense there are the following correspondences:
\dt{24q}{rcl}{
vertex-disjoint directed paths & $\longleftrightarrow$ & free groups, \\
arc-disjoint directed paths & $\longleftrightarrow$ & free abelian groups.
}
In the undirected case we could add the relations $g_i^{2}=1$ for $i=1,\ldots,k$.
This gives the free boolean group (all words made from $g_1,g_2,\ldots$ with no segment $g_ig_i$ for any $i$) and the free abelian boolean groups ($\cong \{0,1\}^k$), and the following correspondences:
\dt{24r}{rcl}{
vertex-disjoint undirected paths & $\longleftrightarrow$ & free boolean groups, \\
edge-disjoint undirected paths & $\longleftrightarrow$ & free abelian boolean groups.
}
By Robertson and Seymour's result, for fixed $k$ the $k$ disjoint undirected paths problem is solvable in polynomial time (for the vertex-disjoint case, and hence also for the edge-disjoint case).
This might suggest that problem \rf{24p} is solvable in polynomial time for any fixed free boolean (abelian) group.

However, problem \rf{24p} is NP-complete for $G:=\{0,1\}^2$ and $C:=\{(0,0),(1,0),(0,1)\}$, even if we fix $\phi(a)=(1,1)$ for each arc $a$.
In that case \rf{24p} has a solution $\psi$ if and only if $D$ is four
vertex colorable.
(I thank Bert Gerards for this observation.)

\vspace{4mm}

\section*{References}\label{REF}
{\small
\begin{itemize}{}{
\setlength{\labelwidth}{4mm}
\setlength{\parsep}{0mm}
\setlength{\itemsep}{1mm}
\setlength{\leftmargin}{5mm}
\setlength{\labelsep}{1mm}
}
\item[\mbox{\rm[1]}] A. Baudisch, 
Kommutationsgleichungen in semifreien Gruppen,
{\em Acta Mathematica A\-ca\-de\-mi\-ae Scientiarum Hungaricae}
29 (1977) 235--249.

\item[\mbox{\rm[2]}] G. Birkhoff, 
{\em Lattice Theory (Third Edition)}
[American Mathematical Colloquium Publications Volume {XXV}],
American Mathematical Society, Providence, Rhode Island, 1973.

\item[\mbox{\rm[3]}] G. Ding, A. Schrijver, P.D. Seymour, 
Disjoint paths in a planar graph --- a general theorem,
{\em {SIAM} Journal on Discrete Mathematics} 5 (1992) 112--116.

\item[\mbox{\rm[4]}] C. Droms, 
Isomorphisms of graph groups,
{\em Proceedings of the American Mathematical Society} 100 (1987) 407--408.

\item[\mbox{\rm[5]}] S. Fortune, J. Hopcroft, J. Wyllie, 
The directed subgraph homeomorphism problem,
{\em Theoretical Computer Science} 10 (1980) 111--121.

\item[\mbox{\rm[6]}] J.F. Lynch, 
The equivalence of theorem proving and the interconnection problem,
{\em ({ACM}) {SIGDA} Newsletter} 5:3 (1975) 31--36.

\item[\mbox{\rm[7]}] R.C. Lyndon, P.E. Schupp, 
{\em Combinatorial Group Theory},
Springer, Berlin, 1977.

\item[\mbox{\rm[8]}] W. Magnus, A. Karrass, D. Solitar, 
{\em Combinatorial Group Theory},
Wiley-Interscience, New York, 1966.

\item[\mbox{\rm[9]}] B.A. Reed, N. Robertson, A. Schrijver, P.D. Seymour, 
Finding disjoint trees in planar graphs in linear time,
in: {\em Graph Structure Theory}
(Proceedings Joint Summer Research Conference on Graph Minors,
Seattle, Washington, 1991; N. Robertson, P. Seymour, eds.)
[Contemporary Mathematics 147],
American Mathematical Society, Providence, Rhode Island, 1993,
pp. 295--301.

\item[\mbox{\rm[10]}] N. Robertson, P.D. Seymour, 
Graph minors. {XIII}. The disjoint paths problem,
{\em Journal of Combinatorial Theory, Series B} 63 (1995) 65--110.

\item[\mbox{\rm[11]}] A. Schrijver, 
Disjoint circuits of prescribed homotopies
in a graph on a compact surface,
{\em Journal of Combinatorial Theory, Series B} 51 (1991) 127--159.

\item[\mbox{\rm[12]}] A. Schrijver, 
Disjoint homotopic paths and trees in a planar graph,
{\em Discrete \& Computational Geometry} 6 (1991) 527--574.

\item[\mbox{\rm[13]}] H. Servatius, 
Automorphisms of graph groups,
{\em Journal of Algebra} 126 (1989) 34--60.

\item[\mbox{\rm[14]}] P.D. Seymour, 
Directed circuits on a torus,
{\em Combinatorica} 11 (1991) 261--273.

\item[\mbox{\rm[15]}] M. Sholander, 
Trees, lattices, and betweenness,
{\em Proceedings of the American Mathematical Society} 3 (1952) 369--381.

\item[\mbox{\rm[16]}] M. Sholander, 
Medians and betweenness,
{\em Proceedings of the American Mathematical Society} 5 (1954) 801--807.

\item[\mbox{\rm[17]}] M. Sholander, 
Medians, lattices and trees,
{\em Proceedings of the American Mathematical Society} 5 (1954) 808--812.

\item[\mbox{\rm[18]}] C. Wrathall, 
The word problem for free partially commutative groups,
{\em Journal of Symbolic Computation} 6 (1988) 99--104.

\end{itemize}
}

\end{document}